\documentclass[12pt]{article}

\usepackage{multirow}
\usepackage{mathtools} 
\usepackage{amsmath, amssymb} 
\usepackage{amsthm}
\usepackage{enumerate}  
\usepackage{url} 

\usepackage{authblk}

\newcommand*\samethanks[1][\value{footnote}]{\footnotemark[#1]}
\usepackage{tikz}
\usepackage{caption}

\usepackage{graphicx}

\newcommand\cx{{\mathbb C}}
\newcommand\fld{{\mathbb F}}
{\tiny }
\newcommand\ints{{\mathbb Z}}
\newcommand\rats{{\mathbb Q}}

\DeclarePairedDelimiter\abs{\lvert}{\rvert}%
\DeclarePairedDelimiter\norm{\lVert}{\rVert}%

\makeatletter
\let\oldabs\abs
\def\abs{\@ifstar{\oldabs}{\oldabs*}}
\let\oldnorm\norm
\def\norm{\@ifstar{\oldnorm}{\oldnorm*}}
\makeatother

\newcommand\comp[1]{{\mkern2mu\overline{\mkern-2mu#1}}}
\newcommand\pmat[1]{\begin{pmatrix} #1 \end{pmatrix}}
\newcommand\seq[4]{#1_{#2},#1_{#3},\ldots,#1_{#4}}

\newtheoremstyle{plainsl}%
	{\topsep}
	{\topsep}
	{\slshape} 
	{}
	{\normalfont\bfseries}
	{.}
	{ }
	{}

\swapnumbers

{\theoremstyle{plainsl}
\newtheorem{theorem}{Theorem}[section]
\newtheorem{lemma}[theorem]{Lemma}
\newtheorem{corollary}[theorem]{Corollary}}
{\theoremstyle{remark}
}

\renewcommand\proof{\noindent\textsl{Proof. }}
\newcommand\sqr[2]{{\vbox{\hrule height.#2pt
    \hbox{\vrule width.#2pt height#1pt \kern#1pt
        \vrule width.#2pt}\hrule height.#2pt}}}
\renewcommand\qed{%
	\ifmmode\eqno\sqr53
	\else\nolinebreak\ \hfill\sqr53\medbreak\fi}

\DeclareMathOperator{\tr}{tr}

\newcommand\one{{\bf1}}

\newcommand\grp[1]{\langle #1\rangle}

\usepackage{blkarray}
\usepackage{pbox}
\usepackage{longtable}

\DeclareMathOperator{\cB}{\mathcal{B}}
\DeclareMathOperator{\cC}{\mathcal{C}}

\title{Discrete-Time Quantum Walks and Graph Structures}
\author{Chris Godsil\thanks{Department of Combinatorics and Optimization, University of Waterloo. \texttt{\{cgodsil, h3zhan\}@uwaterloo.ca}} \quad Hanmeng Zhan \samethanks[1] \thanks{Currently a postdoctoral fellow at the Centre de Recherches Math\'ematiques, Universit\'e de Montr\'eal. \texttt{zhanhanm@crm.umontreal.ca}}}

\begin{document}
\maketitle
\begin{abstract}
We formulate three current models of discrete-time quantum walks in a combinatorial way. These walks are shown to be closely related to rotation systems and 1-factorizations of graphs. For two of the models, we compute the traces and total entropies of the average mixing matrices for some cubic graphs. The trace captures how likely a quantum walk is to revisit the state it started with, and the total entropy measures how close the limiting distribution is to uniform. Our numerical results indicate three relations between quantum walks and graph structures: for the first model, rotation systems with higher genera give lower traces and higher entropies, and for the second model, the symmetric 1-factorizations always give the highest trace. 
\end{abstract}


\section{Introduction}
There are at least three different models of discrete-time quantum walks in current research. Some of them are equivalent to quantum search algorithms \cite{Shenvi2003}, and some are shown to exhibit desired properties such as periodicity and perfect state transfer \cite{Barr2014,Dukes2014}. Most of the existing studies focus on a specific model or a restricted family of graphs. In this paper, we analyze the dependence of these models on different combinatorial structures. To be more specific, we reformulate all three models in graph-theoretic language, and experiment on different classes of graphs to see how the properties of a graph affect the properties of a walk.

Unlike the continuous-time quantum walks, the discrete ones require more than just a graph to build. The extra inputs might include a selection of coins, a set of linear orders of the neighbors of vertices, a 1-factorization of the bipartite double cover, or a Markov chain on the graph. According to how the walks are constructed from these building blocks, we will name them the \textsl{arc-reversal model}, the \textsl{shunt-decomposition model} and the \textsl{two-reflection model}. 

In the arc-reversal model, the walk alternately flips a unitary coin and reverses the arcs of the graph. We show that if the coins are circulants with simple eigenvalues, then for each rotation system of the graph, we can construct a unique arc-reversal walk. A detailed discussion on graph embeddings and rotation systems is given in Section \ref{sec_maps} and Section \ref{sec_rotsys}. In the shunt-decomposition model, the walk alternately flips a unitary coin and hops between outgoing arcs in the same class. Each such class, or ``shunt",  determines a perfect matching of the bipartite double cover of the graph. Thus every shunt-decomposition walk inherits the properties of some 1-factorization. The two-reflection model is coinless, and the walk alternately applies two involutions with respect to two weighted partitions of the arcs. These partitions usually arise from a classical Markov chain.

After introducing different models, we review some asymptotic properties of a general quantum walk. It is well-known that the probability distribution of a quantum walk does not converge. However, the average probabilities over time converge to a limit, which we call the average mixing matrix. In Section \ref{avgmix}, we develop some theory of the average mixing matrix. Two parameters of the average mixing matrix are of interest to us---the trace, which captures how likely a walk is to stay at home in the limit, and the total entropy, which captures how far the limiting distributions are from being uniform. 

We then present some computational results indicating sensitivity of walk properties to details of the model. For each cubic graph on up to 12 vertices, we enumerate all its rotation systems and shunt-decompositions, compute the associated average mixing matrices for the arc-reversal model and the shunt-decomposition model, and compared their traces and total entropies. Our numerical results on the traces indicate two relations. For the first model, as the genus of the embedding increases, the trace of the average mixing matrix decreases, with only few exceptions for some cubic graphs on 12 vertices. For the second model, the trace distinguishes non-isomorphic shunt-decompositions, and the symmetric shunts, if any, always give the highest trace. Meanwhile, for the arc-reversal model, the total entropy distinguishes different graph embeddings, and increases as the genus increases for most of the time.

Besides the natural question arising from the aforementioned experiments, we list some other open problems in the end. It is not our current goal to give a charaterization of the most general form of a model. Rather, we aim to find some specific choice of the building blocks, such as the coins, which induce interesting walks that will reveal clear relation between the properties of the graphs and the properties of the walks. 

\section{Models}
A discrete-time quantum walk, in its most general form, is a power of a unitary matrix, which updates the amplitudes on arcs at each step. While it seems natural to construct a unitary matrix acting on the vertices of a graph, not every graph admits a unitary weighted adjacency matrix. The path on three vertices is an example. In \cite{Severini2003}, Severini discovered some properties of digraphs that admit unitary representations. Alternatively, for an arbitrary graph we can always construct a unitary matrix acting on a larger space. In each of the models we describe below, the transition matrix is acting on the arcs, or the ordered pairs of vertices, and is a product of two sparse unitary matrices. The first two types of walks are driven by a set of coins, while the last one is coinless.

\subsection{Arc-Reversal Model}
The first model we investigate dates back to 2001, when Watrous \cite{Watrous2001} presented a method to simulate random walks on regular graphs using quantum computations. In \cite{Kendon2003} Kendon extended his idea to a class of coined quantum walks on general graphs. Based on this model, Emms et al \cite{Emms2005,Emms2009} proposed an algorithm for the graph isomorphism problem. They tested a large number of strongly regular graphs, and the spectrum of their proposed matrix successfully distinguished all of their non-isomorphic strongly regular graphs with the same parameters. 

We will view an undirected graph $X$ as a directed graph, with each edge replaced by a pair of opposite arcs. Suppose $X$ has $n$ vertices and $m$ arcs. Our quantum walk takes place in the complex inner product space $\cx^{m}$, spanned by the characteristic vectors $e_{u,v}$ of the arcs $(u,v)$. Let $R$ be the permutation matrix that reverses each arc, that is, 
\[R e_{u,v} = e_{v,u}.\]
This gives us the first sparse unitary matrix. To construct a second unitary, for each vertex $u$ we specify a linear order on its neighbors:
\[f_u:\{1,2,\cdots,\deg(u)\} \to \{v: u\sim v\}.\]
The vertex $f_u(j)$ will be referred to as the $j$-th neighbor of $u$, and the arc $(u, f_u(j))$ $j$-th arc of $u$. Next, define a $\deg(u) \times \deg(u)$ unitary matrix $C_u$ that acts on the outgoing arcs of $u$. This matrix serves as a quantum coin: it sends the $j$-th arc of $u$ to a superposition of all outgoing arcs of $u$, in which the amplitudes come from the $j$-th column of $C_u$:
\[C_u e_j = \sum_{k=1}^{\deg(u)} (e_k^T C_u e_j) e_k.\]
Denote by $C$ the block diagonal matrix 
\[C = \pmat{C_1 &&&\\ & C_2 &&\\ &&\ddots &\\ &&& C_n}.\]
The transition matrix of an arc-reversal quantum walk is given by
\[U = RC.\]
Thus, each iteration consists of two steps---a coin flip and an arc-reversal.

By way of example, consider the complete graph $K_3$. Pick the linear orders $f_1$, $f_2$ and $f_3$ such that 
\begin{gather*}
f_1(1)=2, \quad f_2(2)=3;\\
f_2(1)=1, \quad f_2(2)=3;\\
f_3(1)=1, \quad f_3(2)=2.
\end{gather*} 
Then the arc-reversal matrix can be written as
\[R=
\begin{blockarray}{ccccccc}
& (1,2) & (1,3) & (2,1) & (2,3) & (3,1) & (3,2) \\
\begin{block}{c(cccccc)}
(1,2) & 0 & 0 & 1 & 0 & 0 & 0\\
(1,3) & 0 & 0 & 0 & 0 & 1 & 0\\
(2,1) & 1 & 0 & 0 & 0 & 0 & 0\\
(2,3) & 0 & 0 & 0 & 0 & 0 & 1\\
(3,1) & 0 & 1 & 0 & 0 & 0 & 0\\
(3,2) & 0 & 0 & 0 & 1 & 0 & 0\\
\end{block}
\end{blockarray}
\]
Now let us give the same coin to all the vertices, so that the $k$-th arc of $u$ hops to the $j$-th arc of $u$ with amplitude $1/\sqrt{2} (-1)^{jk}$, for any $u$. The coin operator is thus
\[C=
\begin{blockarray}{ccccccc}
& (1,2) & (1,3) & (2,1) & (2,3) & (3,1) & (3,2) \\
\begin{block}{c(cccccc)}
(1,2) & 1/\sqrt{2} & 1/\sqrt{2} & 0 & 0 & 0 & 0\\
(1,3) & 1/\sqrt{2} & -1/\sqrt{2} & 0 & 0 & 0 & 0\\
(2,1) & 0 & 0 & 1/\sqrt{2} & 1/\sqrt{2} & 0 & 0\\
(2,3) & 0 & 0 & 1/\sqrt{2} & -1/\sqrt{2} & 0 & 0\\
(3,1) & 0 & 0 & 0 & 0 & 1/\sqrt{2} & 1/\sqrt{2}\\
(3,2) & 0 & 0 & 0 & 0 & 1/\sqrt{2} & -1/\sqrt{2}\\
\end{block}
\end{blockarray}\]
With initial state $e_{1,3}$, after one iteration, the system will be in the state
\[\frac{1}{\sqrt{2}} e_{2,1} - \frac{1}{\sqrt{2}} e_{3,1}.\]

The above coin is called the \textsl{Fourier coin}:
\[F= \left(\frac{1}{\sqrt{d}}e^{2jk\pi i/d}\right)_{jk}\]
Another common choice would be the \textsl{Grover coin}:
\[G = \frac{2}{d}J - I.\] 
This coin treats the neighbors in a simple way -- an arc $(u,v)$ stays at current position with amplitude $2/d$, and moves to any other outgoing arc of $u$ with amplitude $2/d-1$. Hence it is indifferent to the linear order $f_u$. The famous algorithm, Grover's search, can be formulated as an arc-reversal walk on the looped $K_n$, with coin $-G$ assigned to the marked vertex, and coin $G$ to the unmarked ones. This reformulation was first made by Szegedy \cite{Szegedy2004}.

\subsection{Shunt-Decomposition Model}
The next model has been mostly applied to graphs with some symmetry, such as the infinite paths, infinite grids, cycles, and cubes. It is formally defined by Aharonov et al \cite{Aharonov2000}. Similar to an arc-reversal model, each iteration is a coin flip followed by another unitary operation. Let $\{f_u: u\in V(X)\}$ be a set of linear orders as described in the previous section. For reasons that will become clear later, assume our graph $X$ is $d$-regular. The state space $\cx^{m}$ is then isomorphic to the tensor product $\cx^n \otimes \cx^d$. With this decomposition, the arc $(u, f_u(j))$ can be represented by the vector
\[ e_u \otimes e_j.\]
The transition matrix is given by
\[U=SC,\]
where $C$ is the coin operator, and $S$ shifts the arcs in the following way: for each vertex $u$, it maps its $j$-th arc to the $j$-th arc of $f_u(j)$: 
\[S ( e_u \otimes e_j) = e_{f_u(j)} \otimes e_j.\]
That is, $S$ moves every arc one step forward to the arc with the same label. By definition, $S$ is a $01$-matrix where each column has exactly one $1$. Since $S$ is unitary, it has to be a permutation matrix. Hence the linear orders $f_u$ cannot be arbitrary; for any vertex $u$, the values
\[\{f_w^{-1}(v): w \sim v\}\]
must be distinct. Under this assumption, each label $j$ induces a permutation on $V(X)$ that maps one vertex to its $j$-th neighbor, denoted by the permutation matrix $P_j$. This explains why $X$ has to be regular.  Now, $S$ is equivalent to the block diagonal matrix
\[\pmat{P_1 & & &\\
	& P_2 & & \\
	& & \ddots &\\
	& & & P_d}.\]
In what follows, we will let $E_{ij}$ denote the matrix with $1$ in the $ij$-entry and $0$ elsewhere. By the isomorphism between $\cx^n\otimes \cx^d$ and $\cx^d \otimes \cx^n$, 
\[U = (P_1\otimes E_{11}+\cdots+P_d\otimes E_{dd}) (E_{11}\otimes C_1+\cdots+E_{nn}\otimes C_n).\]

We consider an example. The following is a collection of linear orders for $K_3$:
\begin{gather*}
f_1(1)=2, \quad f_1(2)=3;\\
f_2(1)=3, \quad f_2(2)=1;\\
f_3(1)=1, \quad f_3(2)=2.
\end{gather*} 
They give rise to an ordered coloring of the arcs of $K_3$, as shown in Figure \ref{arc_coloring}. 

\begin{figure}[!htb]
\centering
\begin{minipage}[b]{0.4\textwidth}
	\centering
	\begin{tikzpicture}
	[every node/.style={circle, draw}]
	
	\node[label=above:$1$] (1) at (0,1.6){};
	\node[label=left:$2$] (2) at (-1,0){};
	\node[label=right:$3$] (3) at (1,0){};
	
	\foreach \a/\b in {1/2, 2/3, 3/1}
	\draw[bend right=18, ->, red] (\a) to (\b);
	
	\foreach \a/\b in {1/3, 3/2, 2/1}
	\draw [bend right=18, ->, blue] (\a) to (\b);
	\end{tikzpicture}	
	\captionof{figure}{an arc coloring of $K_3$}
	\label{arc_coloring}
\end{minipage}%
\begin{minipage}[b]{0.6\textwidth}
	\centering
	\begin{tikzpicture}
	[every node/.style={circle, draw}]
	
	\node[label=left:$1$] (1) at (0,1){};
	\node[label=left:$2$] (2) at (0,0){};
	\node[label=left:$3$] (3) at (0,-1){};
	
	\node[label=right:$1'$] (1') at (2,1){};
	\node[label=right:$2'$] (2') at (2,0){};
	\node[label=right:$3'$] (3') at (2,-1){};
	
	\foreach \a/\b in {1/2', 2/3', 3/1'}
	\draw[red] (\a) -- (\b);
	
	\foreach \a/\b in {1/3', 2/1', 3/2'}
	\draw[blue] (\a) -- (\b);
	\end{tikzpicture}
	\captionof{figure}{an edge-coloring of $K_2 \times K_3$}
	\label{1-factn}
\end{minipage}
\end{figure}
Each color determines a permutation on the vertices:
\[P_1 =\begin{blockarray}{cccc}
& 1 & 2 & 3  \\
\begin{block}{c(ccc)}
1 & 0 & 0 & 1\\
2 & 1 & 0 & 0\\
3 & 0 & 1 & 0\\
\end{block}
\end{blockarray}
\qquad P_2 =\begin{blockarray}{cccc}
& 1 & 2 & 3  \\
\begin{block}{c(ccc)}
1 & 0 & 1 & 0\\
2 & 0 & 0 & 1\\
3 & 1 & 0 & 0\\
\end{block}
\end{blockarray}\]
which sum up to the adjacency matrix of $K_3$. Equivalently, we can realize the linear orders as an ordered edge coloring of the bipartite graph $K_2 \times K_3$, as shown in Figure \ref{1-factn}. The edge $(j,k')$ is colored red if $f_j(1)=k$, and blue if $f_j(2)=k$.

In general, for a $d$-regular graph $X$, there is a one-to-one correspondence between a decomposition of its adjacency matrix into $d$ permutation matrices
\[A(X) = P_1 + \cdots + P_d\]
and a 1-factorization of $K_2\times X$. Each permutation $P_j$ maps vertices to adjacent vertices. Such permutations are called \textsl{shunts}. We refer to the above decomposition as a \textsl{shunt-decomposition} of $X$. A set of coins together with an ordered shunt-decomposition determines a quantum walk.

\subsection{Two-Reflections Model}

The third model is a direct quantization of the classical random walks, without introducing a coin register. It was extracted from Ambainis's paper \cite{Ambainis2003} and formalized by Szegedy \cite{Szegedy2004}. Let $M$ be a Markov chain on $X$. We will construct two reflections $R_1$ and $R_2$ on the $n^2$ ordered pairs of vertices based on $M$, and set 
\[U = R_2R_1.\]
First, let $N$ denote the matrix obtained from $M$ by taking the square root of its entries. Next, define two partitions of the ordered pairs:
\begin{align*}
Q_1 &= \pmat{e_1\otimes (Ne_1) & \cdots & e_n\otimes (Ne_n)}\\
Q_2 &= \pmat{(N^Te_1) \otimes e_1 & \cdots & (N^Te_n) \otimes e_n}.
\end{align*}
Since $M$ is doubly stochastic, these are $n^2\times n$ matrices such that 
\[Q_1^TQ_1 = Q_2^TQ_2 = I.\]
Then for each $j$, the matrix $Q_jQ_j^T$ is a projection onto the column space of $Q_j$, and 
\[R_j = 2Q_jQ_j^T - I\]
is a reflection about the column space of $Q_j$.

For instance, take the simple random walk on $K_3$ with
\[M = \begin{blockarray}{cccc}
& 1 & 2 & 3  \\
\begin{block}{c(ccc)}
1 & 0 & 1/2 & 1/2\\
2 & 1/2 & 0 & 1/2\\
3 & 1/2 & 1/2 & 0\\
\end{block}
\end{blockarray}.\]
Then we have
\[Q_1 = \begin{blockarray}{cccc}
& 1 & 2 & 3  \\
\begin{block}{c(ccc)}
11 & 0 & 0 & 0\\
12 & 1/\sqrt{2} & 0 & 0\\
13 & 1/\sqrt{2} & 0 & 0\\
21 & 0 & 1/\sqrt{2} & 0\\
22 & 0 & 0 & 0\\
23 & 0 & 1/\sqrt{2} & 0\\
31 & 0 & 0 & 1/\sqrt{2}\\
32 & 0 & 0 & 1/\sqrt{2}\\
33 & 0 & 0 & 0\\
\end{block}
\end{blockarray}
\qquad
Q_2 = \begin{blockarray}{cccc}
& 1 & 2 & 3  \\
\begin{block}{c(ccc)}
11 & 0 & 0 & 0\\
12 & 0 & 1/\sqrt{2} & 0\\
13 & 0 & 0 & 1/\sqrt{2}\\
21 & 1/\sqrt{2} & 0 & 0\\
22 & 0 & 0 & 0\\
23 & 0 & 0 & 1/\sqrt{2}\\
31 & 1/\sqrt{2} & 0 & 0\\
32 & 0 & 1/\sqrt{2} & 0\\
33 & 0 & 0 & 0\\
\end{block}
\end{blockarray}
\]
Effectively, $Q_1$ partitions the arcs according to their tails, and $Q_2$ partitions the arcs according to their heads. This is true for any simple random walk, as for any non-adjacent vertices $u$ and $v$, the row indexed by $(u,v)$ is zero in both $Q_1$ and $Q_2$. Thus, a quantization of a simple random walk can be viewed as a walk on the arcs only. In the above example, the initial state $e_{12}$ is mapped to $e_{13}$ by $R_1$, and then mapped to $e_{12}$ by $R_2$.

More generally, we can define a two-reflection walk by setting 
\[U = (2Q_1Q_1^T - I) (2Q_2 Q_2^T - I),\]
where $Q_1$ and $Q_2$ are normalized character matrices that represent two partitions of the arcs. These partitions may arise from other graph structures, such as orientable embeddings \cite{Zhan2017a}. The spectral decomposition of a two-reflection walk has been studied \cite{Szegedy2004}, and used to find perfect state transfer in circulant graphs \cite{Zhan2017}.

\subsection{Similarities and Differences \label{subsec_simdiff}}
In some cases the aforementioned models are equivalent, but there are also quantum walks that belong to only one of the frameworks. The shift matrix $S$ in the shunt-decomposition is symmetric if and only if every shunt $P_j$ is symmetric, and so the arc coloring arising from the shunt-decomposition is in fact an edge coloring. Hence, the arc-reversal model and the shunt-decomposition model overlap when the graph admits a $d$-edge-coloring, and the same coins and same linear orders are applied. The arc-reversal model can also be equivalent to the two-reflections model, if the first walk uses the Grover coins and lexicographical linear orders, and the second walk quantizes the simple random walk. They are equivalent in the sense that two steps of the first walk has the same effect as one step of the second walk. For more details, see the recent paper by Wong \cite{Wong2016}. Another correlation between these two models has been found by Portugal and Segawa \cite{Portugal2016}, who showed that an arc-reversal walk with Grover coins can be converted into a two-reflection walk on the same graph with each edge subdivided once.

\section{Shunt-Decompositions}
We have seen that a shunt-decomposition of a graph $X$ is equivalent to a 1-factorization of $K_2 \times X$. Every regular bipartite graph has at least one 1-factorization, but there might be more, in which case we would like to compare the properties of walks induced by different shunt-decompositions. To enumerate all 1-factorizations of a graph $Y$, notice that a matching of $Y$ is a coclique of the line graph $L(Y)$ of $Y$. Thus, every 1-factorization of $Y$ gives rise to a vertex partition of $L(Y)$ into maximum cocliques, and vise versa. We enumerated all shunt-decompositions of each cubic graph on up to 8 vertices, and compared the limiting average distributions of the corresponding quantum walks. More details will be discussed in Section \ref{data_shunt}.

Some graphs come with natural shunt-decompositions. For example, the arcs in a Cayley graph $X(G, \cC)$ are partitioned according to the elements in the connection set $\cC$. Hence each element $c_j \in \cC$ determines a shunt that sends vertex $u$ to vertex $c_i^{-1} u$. In the next section, we will consider another family of graphs with obvious shunt-decompositions.

\section{Maps \label{sec_maps}}
Let $X$ be a graph embedded on some surface. For simplicity, we assume that every face is bounded by a cycle. Let $u$ be a vertex, $e$ an edge incident to $u$, and $f$ a face containing the edge $e$. The triple $(u,e,f)$ is called a flag of the embedding. Pictorially, a flag is a triangle in the barycentric subdivision of a face. In the following example,  each black dot represents a flag in the planar embedding of $C_3$.

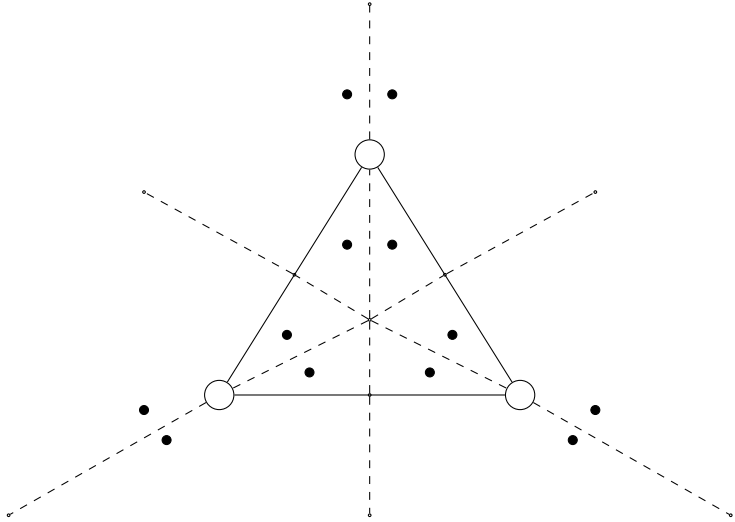
\begin{figure}[!htb]
\centering
\begin{tikzpicture}
[every node/.style={circle,draw}]

\node (1) at (0,3.2){};
\node (2) at (-2,0){};
\node (3) at (2,0){};

\tikzset{inv/.style={inner sep=0pt}}

\node[inv] (0) at (0,1){};
\node[inv] (1') at (0,0){};
\node[inv] (2') at (1,1.6){};
\node[inv] (3') at (-1,1.6){};
\node[inv] (a) at (0,5.2){};
\node[inv] (a') at (0,-1.6){};
\node[inv] (b) at (-4.8,-1.6){};
\node[inv] (b') at (3,2.7){};
\node[inv] (c) at (4.8,-1.6){};
\node[inv] (c') at (-3,2.7){};

\tikzset{dot/.style={draw, inner sep=1.2pt, fill}}

\node[dot] (A) at (-0.3,2){};
\node[dot] (B) at (0.3,2){};
\node[dot] (C) at (1.1,0.8){};
\node[dot] (D) at (0.8,0.3){};
\node[dot] (E) at (-0.8,0.3){};
\node[dot] (F) at (-1.1,0.8){};

\node[dot] (G) at (-0.3,4){};
\node[dot] (H) at (0.3,4){};
\node[dot] (I) at (3,-0.2){};
\node[dot] (J) at (2.7,-0.6){};
\node[dot] (K) at (-2.7,-0.6){};
\node[dot] (L) at (-3,-0.2){};

\foreach \a/\b in {1/2, 2/3, 3/1}
\draw (\a) to (\b);

\foreach \a/\b in {0/1, 0/1', 0/2, 0/2', 0/3, 0/3', 1/a,1'/a',2'/b', 3'/c',2/b,3/c}
\draw [dashed] (\a) to (\b);

\end{tikzpicture}
\caption{Flags of $C_3$ are represented by the black dots}
\end{figure}

For each flag $(u,e,f)$, let $u'$ be the other endpoint of $e$, let $e'$ be the other edge in $f$ that is incident to $u$, and let $f'$ be the other face that contains $e$. Define three functions
\begin{align*}
\tau_0: &(u, e, f) \mapsto (u', e, f),\\
\tau_1: &(u, e, f) \mapsto (u, e', f),\\
\tau_2: &(u, e, f) \mapsto (u, e, f').
\end{align*}
We have the following observations.
\begin{enumerate}[(i)]
\item $\tau_0, \tau_1, \tau_2$ are fixed-point-free involutions.
\item $\tau_0\tau_2 = \tau_2\tau_0$, and $\tau_0\tau_2 $ is fixed-point-free.
\item The group $\grp{\tau_0, \tau_1, \tau_2}$ acts transitively on the flags.
\end{enumerate}
Combinatorially, any tuple $(\tau_0, \tau_1, \tau_2)$ satisfying the above is called a map. More details on maps can be found in \cite{Archdeacon1994}.

\begin{figure}[!htb]
\centering
\begin{tikzpicture}
[every node/.style={circle,draw}]

\node (1) at (0,3.2){};
\node (2) at (-2,0){};
\node (3) at (2,0){};

\tikzset{inv/.style={inner sep=0pt}}

\node[inv] (0) at (0,1){};
\node[inv] (1') at (0,0){};
\node[inv] (2') at (1,1.6){};
\node[inv] (3') at (-1,1.6){};
\node[inv] (a) at (0,5.2){};
\node[inv] (a') at (0,-1.6){};
\node[inv] (b) at (-4.8,-1.6){};
\node[inv] (b') at (3,2.7){};
\node[inv] (c) at (4.8,-1.6){};
\node[inv] (c') at (-3,2.7){};

\tikzset{dot/.style={draw, inner sep=1.2pt, fill}}

\node[dot] (A) at (-0.3,2){};
\node[dot] (B) at (0.3,2){};
\node[dot] (C) at (1.1,0.8){};
\node[dot] (D) at (0.8,0.3){};
\node[dot] (E) at (-0.8,0.3){};
\node[dot] (F) at (-1.1,0.8){};

\node[dot] (G) at (-0.3,4){};
\node[dot] (H) at (0.3,4){};
\node[dot] (I) at (3,-0.2){};
\node[dot] (J) at (2.7,-0.6){};
\node[dot] (K) at (-2.7,-0.6){};
\node[dot] (L) at (-3,-0.2){};

\foreach \a/\b in {1/2, 2/3, 3/1}
\draw (\a) to (\b);

\foreach \a/\b in {0/1, 0/1', 0/2, 0/2', 0/3, 0/3', 1/a,1'/a',2'/b', 3'/c',2/b,3/c}
\draw [dashed] (\a) to (\b);

\foreach \a/\b in{B/C, D/E, F/A, H/I, J/K, L/G}
\draw[red, thick] (\a) to (\b);

\foreach \a/\b in {A/B, C/D, E/F, G/H, I/J, K/L}
\draw[blue, thick] (\a) to (\b);

\foreach \a/\b in {A/G, B/H, C/I, D/J, E/K, F/L}
\draw[green!70!black, thick] (\a) to (\b);
\end{tikzpicture}
\caption{Planar embedding of $C_3$ and its gem}
\label{flag_C3}
\end{figure}
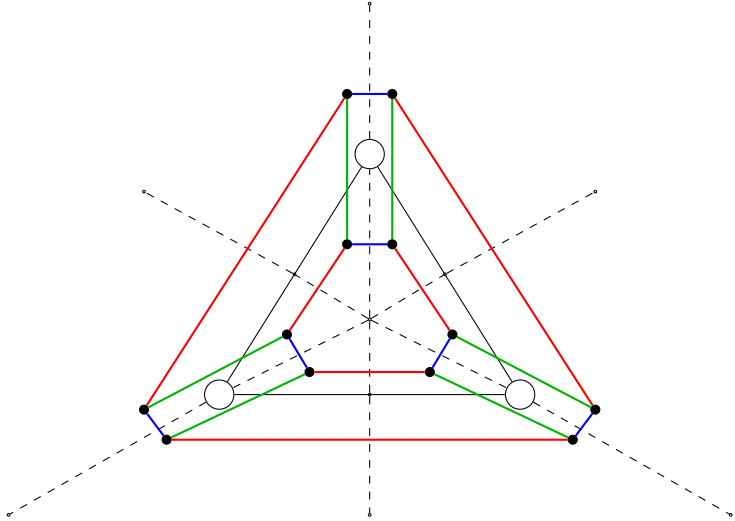

Coming back to our map from an embedding of $X$, the \textsl{graph-encoded-map}, or \textsl{gem}, is the graph whose vertices are the flags, and two flags are adjacent if they are swapped by one of $\tau_0, \tau_1, \tau_2$. We will denote this graph by $\Phi(X)$. Note that $\Phi(X)$ is cubic and $3$-edge-colorable. In Figure \ref{flag_C3}, we draw the gem of the planar $C_3$ with red, blue and green edges. The gem of the planar $C_n$ is the prism graph $K_2\square C_{2n}$.

\begin{figure}[!htb]
	\centering
	\begin{tikzpicture}
	[every node/.style={circle}]
	
	\node[draw] (0) at (0,0) {};
	\node[draw] (1) at (2,0) {};
	\node[draw] (2) at (2,2) {};
	\node[draw] (3) at (0,2) {};
	\node[draw] (4) at (4,0) {};
	\node[draw] (5) at (4,2) {};
	
	\node (0rt) at (0.2,0.2) {};
	\node (0lt) at (-0.2,0.2) {};
	\node (0rb) at (0.2,-0.2) {};
	\node (1rt) at (2.2,0.2) {};
	\node (1lt) at (1.8,0.2) {};
	\node (1rb) at (2.2,-0.2) {};
	\node (1lb) at (1.8,-0.2) {};
	\node (2rt) at (2.2,2.2) {};
	\node (2lt) at (1.8,2.2) {};
	\node (2rb) at (2.2,1.8) {};
	\node (2lb) at (1.8,1.8) {};
	\node (3rt) at (0.2,2.2) {};
	\node (3rb) at (0.2,1.8) {};
	\node (3lb) at (-0.2,1.8) {};
	\node (4rt) at (4.2,0.2) {};
	\node (4lt) at (3.8, 0.2) {};
	\node (4lb) at (3.8,-0.2) {};
	\node (5lt) at (3.8,2.2) {};
	\node (5rb) at (4.2,1.8) {};
	\node (5lb) at (3.8,1.8) {};
	
	\foreach \a/\b in {0/1, 1/2, 2/3, 3/0, 1/4, 4/5, 5/2}
	\draw (\a) to (\b);
	
	\foreach \a/\b in {0rt/3rb, 3rb/2lb,2lb/1lt, 1lt/0rt, 2rb/5lb, 5lb/4lt, 4lt/1rt, 1rt/2rb, 3lb/0lt, 0rb/1lb, 1rb/4lb, 4rt/5rb, 5lt/2rt, 2lt/3rt}
	\draw[->] (\a) to (\b);
	\end{tikzpicture}
	\caption{Coherent orientation of faces}
	\label{coherent}
\end{figure}
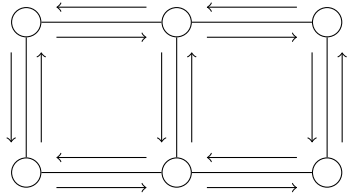

The map $(\tau_0, \tau_1, \tau_2)$  is \textsl{orientable} if the underlying surface is orientable. For a formal definition of orientable manifolds, see Lee \cite[Ch 15]{Lee2000}. Here we only provide some intuition. Given a face on a surface, we can orient its boundary in two directions. If there is some orientation of the faces such that whenever an edge $e$ is contained in two faces $f$ and $f'$, the orientation $e$ receives in $f$ is opposite to the orientation it receives in $f'$, then the orietation is called \textsl{coherent}. Figure \ref{coherent} is an example of a coherent orietation. A surface is called \textsl{orientable} if it admits a coherent orientation, and \textsl{non-orientable} otherwise.

 Alternatively, we can characterize orientability using the gem. The following result is due to Vince \cite{Vince1983}.
 
 \begin{theorem}[Vince \cite{Vince1983}]
 The embedding of $X$ is orientable if and only if the gem $\Phi(X)$ is bipartite.
 \qed
 \label{orient_bip}
 \end{theorem}

With the $3$-edge-coloring induced by $\tau_0, \tau_1, \tau_2$, we obtain a symmetric shunt-decomposition of $\Phi(X)$:
\[A(\Phi(X)) = P_0 + P_1 + P_2,\]
where $P_j^T = P_j$. According to Subsection \ref{subsec_simdiff}, this defines an arc-reversal walk on $\Phi(X)$. Suppose $X$ is embedded on an orientable surface. Then $\Phi(X)$ is bipartite and cubic. If, in addition, $\Phi(X)$ is isomorphic to $K_2\times Y$ for some cubic graph $Y$, then the edge-coloring of $\Phi(X)$ also gives rise to a shunt-decomposition, possibly asymmetric, on the smaller graph $Y$. To illustrate this, recall that the gem of the planer $C_3$ is $K_2\square C_6$. It is isomorphic to $K_2\times (K_2\square K_3)$, as drawn in Figure \ref{flagrf}. Figure \ref{K2_K_3} shows the corresponding shunt-decomposition of the graph $K_2\square K_3$.

\begin{figure}[!htb]
\centering
\begin{minipage}[b]{0.5\textwidth}
	\centering
	\begin{tikzpicture}
	[every node/.style={circle, draw}]
	\node[label=left:$1$] (1) at (-0.9,1.9){};
	\node[label=right:$3'$] (3') at (0.9,1.9){};
	\node[label=right:$2$] (2) at (2.1,0){};
	\node[label=left:$2'$] (2') at (-2.1,0){}; 
	\node[label=left:$3$] (3) at (-0.9,-1.9){};
	\node[label=right:$1'$] (1') at (0.9,-1.9){};
	\node[label={[xshift=0.5em, yshift=-2.45em]$6'$}] (6') at (-0.5,0.9){};
	\node[label={[xshift=-0.5em,yshift=-2.4em]$5$}] (5) at (0.5,0.9){};
	\node[label={[xshift=0.9em,yshift=-1.5em]$4$}] (4) at (-1.1,0){};
	\node[label={[xshift=-0.9em,yshift=-1.5em]$4'$}] (4') at (1.1,0){};
	\node[label={[xshift=0.5em,yshift=-0.7em]$5'$}] (5') at (-0.5,-0.9){};
	\node[label={[xshift=-0.4em,yshift=-0.6em]$6$}] (6) at (0.5,-0.9){};
	
	\foreach \a/\b in {1/2', 2/3',3/1', 4/6',6/5',5/4'}
	\draw[red] (\a) to (\b);
	\foreach \a/\b in {2/1',3/2',1/3',6/4',5/6',4/5'}
	\draw[blue] (\a) to (\b);
	\foreach \a/\b in {1/6', 6/1', 2/4', 4/2', 3/5', 5/3'}
	\draw[green!70!black] (\a) to (\b);
	\end{tikzpicture}
	\captionof{figure}{$K_2\square C_6 \simeq K_2\times (K_2\square K_3)$}
	\label{flagrf}
\end{minipage}%
\begin{minipage}[b]{0.5\textwidth}
	\centering
	\begin{tikzpicture}
	[every node/.style={circle, draw}]
	
	\node[label=above:$1$] (1) at (0,1.9){};
	\node[label=left:$2$] (2) at (-1,0.3){};
	\node[label=right:$3$] (3) at (1,0.3){};
	\node[label=left:$4$] (4) at (-1,-1.9){};
	\node[label=right:$5$] (5) at (1,-1.9){};
	\node[label=below:$6$] (6) at (0,-0.3){};
	
	\foreach \a/\b in {1/2, 2/3, 3/1, 4/6, 6/5, 5/4}
	\draw [bend right=18, ->, red] (\a) to (\b);
	\foreach \a/\b in {2/1, 3/2, 1/3, 6/4, 4/5, 5/6}
	\draw [bend right=18, ->, blue] (\a) to (\b);
	\foreach \a/\b in {1/6, 2/4, 3/5, 6/1, 4/2, 5/3}
	\draw[green!70!black] (\a) to (\b);
	\end{tikzpicture}	
	\captionof{figure}{$K_2\square K_3$}
	\label{K2_K_3}
\end{minipage}
\end{figure}

In general, the correspondence between $\Phi(X)$ and $Y$ is not unique. However, for each $Y$ such that $\Phi(X)$ is isomorphic to $K_2\times Y$, the adjacency matrix of $\Phi(X)$ can be written as 
\[A(\Phi(X))= \pmat {0 & A(Y) \\ A(Y) & 0}.\]
Given the shunt decomposition of $A$:
\[A(\Phi(X))=P_0+P_1+P_2\]
where each $P_j$ is of the form
\[P_j = \pmat{0 & Q_j^T \\ Q_j & 0},\]
we have the shunt-decomposition of $Y$:
\[A(Y) = Q_0+Q_1+Q_2 = Q_0^T+Q_1^T+Q_2^T.\]
Hence, the quantum walk on $Y$ is 
\[U(Y) = (P_0\otimes E_{00} + P_1\otimes E_{11} + P_2\otimes E_{22}) C_Y,\]
while the quantum walk on $\Phi(X)$is
\[U(\Phi(X)) = (Q_0\otimes E_{00} + Q_1\otimes E_{11}+Q_2\otimes E_{22}) C_{\Phi(X)}.\]
With a proper choice of the coins, we will have
\[U(\Phi(X)) = \pmat{0 & S_Y C_Y\\ S_Y^T C_Y & 0}.\]

\section{Rotation Systems \label{sec_rotsys}}
In either of the coined models, we need a set of linear orders $\{f_u: u\in V(X)\}$ of the neighbors to implement the coin operator. If we convert these linear orders into cyclic permutations,  we obtain a rotation system, which determines an orientable embedding of $X$. We introduce some basic concepts about rotation systems. For more background, see Gross and Tucker \cite{Gross2001}.

Formally, a \textsl{rotation system} is a set $\{\pi_u: u\in V(X)\}$ where each $\pi_u$ is a cyclic permutation on the neighbors of the vertex $u$. For any arc $(u_1, u_2)$, consider the walk
\[(u_1, u_2), (u_2, u_3), (u_3, u_4), \cdots, (u_{k-1}, u_k), \cdots\]
where 
\[u_{j+1} = \pi_{u_j} (u_{j-1}).\]
Since the graph is finite, eventually the walk will meet an arc that is already taken. Moreover, the first arc that is used twice must be $(u_1, u_2)$, as the preimage $\pi_u^{-1} (v)$ is uniquely determined for each $u$. We will call the closed walk
\[(u_1, u_2), (u_2, u_3), (u_3, u_4),\cdots, (u_k, u_1)\]
a \textsl{facial walk}. The facial walks partition the arcs of $X$. Clearly, each edge occurs either once in two facial walks, or twice in the same facial walk. On the other hand, a facial walk can use a vertex many times without using an arc twice. The following three embeddings illustrate the difference. In Figure \ref{face_cyc}, every face is bounded by a cycle. In Figure \ref{face_edge}, the outer face uses the pendent edge twice. In Figure \ref{face_vx}, the outer face uses the central vertex four times but no edge twice.

\begin{figure}[!htb]
\centering
\begin{minipage}[b]{0.3\textwidth}
	\centering
	\begin{tikzpicture}
	[every node/.style={circle, draw}]
	
	\node (1) at (1,0){};
	\node (2) at (0,-1){};
	\node (3) at (0,1){};
	\node (4) at (-1,0){};
	
	\foreach \a/\b in {1/2, 2/3, 1/3, 2/4, 3/4}
	\draw (\a) to (\b);
	\end{tikzpicture}	
	\captionof{figure}{every face is bounded by a cycle}\label{face_cyc}
\end{minipage}%
\begin{minipage}[b]{0.3\textwidth}
	\centering
	\begin{tikzpicture}
	[every node/.style={circle, draw}]
	
	\node (1) at (0,1){};
	\node (2) at (0,0){};
	\node (3) at (-1,-1){};
	\node (4) at (1,-1){};
	
	\foreach \a/\b in {1/2, 2/3, 3/4, 2/4}
	\draw (\a) to (\b);
	\end{tikzpicture}
	\captionof{figure}{one face uses an edge twice}\label{face_edge}
\end{minipage}%
\begin{minipage}[b]{0.4\textwidth}
	\centering
	\begin{tikzpicture}
	[every node/.style={circle, draw}]
	
	\node (1) at (0,0){};
	\node (2) at (-1,1){};
	\node (3) at (-1,-1){};
	\node (4) at (1,1){};
	\node (5) at (1,-1){};
	
	\foreach \a/\b in {1/2,1/3,1/4,1/5,2/3,4/5}
	\draw (\a) to (\b);
	\end{tikzpicture}
	\captionof{figure}{one face uses a vertex four times}\label{face_vx}
	\end{minipage}
\end{figure}

Given a facial walk of length $k$, we associate it with a polygon with $k$ sides, labeled by the arcs in the same order as they appear in the facial walk. We then ``glue'' every two sides of these polygons labeled by the same edge. The result is an embedding of the graph onto an orientable surface.

Now let $X$ be a $d$-regular graph and try to construct a quantum walk from a rotation system of $X$. For simplicity, we will choose the arc-reversal model. The first problem is that there are $d$ ways to convert a cyclic permutation into a linear order, so we wish to unify the coin matrix with respect to all possible conversions. To be more specific, if $f_u$ and $g_u$ are two linear orders obtained from the cyclic order $\pi_u$, then the amplitude from arc $(u, f_u(j))$ to arc $(u, f_u(k))$ is required to equal the amplitude from arc $(u, g_u(j))$ to arc $(u, g_u(k))$.

\begin{lemma}
Let $X$ be a $d$-regular graph. Let $f_u$ and $g_u$ be two linear orders of the neighbors of $u$. Let $C_u$ be a unitary coin indexed by the outgoing arcs of $u$. The following are equivalent.
\begin{enumerate}[(i)]
\item For $j,k=1,2,\cdots,d$, the amplitude from arc $(u, f_u(j))$ to arc $(u, f_u(k))$ is equal to the amplitude from arc $(u, g_u(j))$ to arc $(u, g_u(k))$.
\item $C_u$ commutes with the permutation matrix
\[P_u = \pmat{e_{g_uf_u^{-1}(1)} & e_{g_uf_u^{-1}(2)} & \cdots e_{g_uf_u^{-1}(d)}}.\]
\end{enumerate}
\end{lemma}
\proof
For notational ease, let us replace $C_u$, $f_u$ and $g_u$ by $C$, $f$ and $g$, respectively. Statement (i) holds if and only if for each $j$ and $k$, 
\[e^T_{f(k)} C e_{f(j)} = e^T_{g(k)} C e_{g(j)}.\]
On the other hand, since $f$ and $g$ are bijections, statement (ii) holds if and only if for each $j$ and $k$, 
\[e^T_{g(k)} (CP) e_{f(j)} = e^T_{g(k)} (PC) e_{f(j)}.\]
It follows from
\[P e_{f(j)} = e_{g(j)}\]
and
\[P^T e_{g(k)} = P^{-1} e_{g(k)} = e_{f(j)}\]
that the above two statements are equivalent.
\qed

Now, if the linear orders $f_u$ and $g_u$ arise from the same cyclic order $\pi_u$, then they differ in only a cyclic permutation. An immediate consequence is that each coin $C_u$ must be cyclic so that the quantum walk for a rotation system is well-defined. This answers our first question. In addition, for  different rotation systems, we wish to obtain different transition matrices, so the coin $C_u$ should not commute with any non-cyclic permutation. We show that this is guaranteed as long as $C_u$ has simple eigenvalues.

\begin{lemma}
Let $C_u$ be a circulant unitary matrix with simple eigenvalues. Then $C_u$ commutes with a permutation matrix $P_u$ if and only if $P_u$ is cyclic.
\end{lemma}
\proof
We use the fact that for a matrix $C$ with simple eigenvalues, its commutant consists of precisely the polynomials in $C$ (see, for example, Horn and Johnson \cite[Theorem 3.2.4.2]{Horn2012}). The result follows since any polynomial of a circulant matrix is also circulant.
\qed

In Section \ref{data_arcrev}, we will pick a $3\times 3$ circulant coin with simple eigenvalues, and investigate the arc-reversal quantum walks for all rotation systems of cubic graphs on up to 10 vertices.

\section{Hitting Times}
The hitting time measures how fast a quantum walk reaches a specific state. Let $x$ be the initial state, and $y$ the target state. A typical choice for $x$ and $y$ would be the characteristic vectors of two arcs in a graph, respectively. Since the state of the system is in a superposition, it is not on an arc until we perform a measurement. However, if we do a complete measurement after each step, the quantum walk will collapse to a classical random walk. Thus, we need a way to define ``reaching" without killing the quantum coherence.

One choice is to let the quantum walk evolve and do a complete measurement only once. The \textsl{one-shot hitting time}, given by Kempe \cite{Kempe2002}, is the time $k$ when a complete measurement returns the target state with high probability:
\[\abs{y^*U^kx}^2\ge 1-\epsilon.\]
For this definition we need to know when to measure the walk. If the chosen probability $1-\epsilon$ is too high, the hitting time may not exist.

Alternatively, we can do a partial measurement with respect to $y$ after each transition, and determine whether the walk hits $y$ and stops, or is in a state orthogonal to $y$ and continues. The partial measurement consists of the projection $yy^*$ onto the space spanned by $y$, and the projection $I-yy^*$ onto the space orthogonal to $y$. Two notions of hitting time based on this measured walk are given in Kempe \cite{Kempe2002} and Krovi et al \cite{Krovi2006}. 
The \textsl{concurrent hitting time}, defined in \cite{Kempe2002}, is the earliest time before which the walks stops with a fixed high probability $1-\epsilon$:
\[h_{\epsilon} = \min\left\{K: \sum_{k=1}^{K} \abs{ y^* U \left((I-yy^*)U\right)^{k-1} x}^2 \ge 1-\epsilon \right\}.\]
Or, we can follow \cite{Krovi2006} and consider the \textsl{expected hitting time}, the average number of steps the walk takes to reach $x$:
\[h = \sum_{k=1}^{\infty} k\abs{ y^* U \left((I-yy^*)U\right)^{k-1} x}^2.\]

\section{Average Limiting Distributions}

In a classical random walk, the probability distribution converges to a stationary distribution under only mild conditions. In a quantum walk, however, the unitary matrix $U$ preserves the difference between the states at two consecutive steps:
\[\norm{Ux-x} = \norm{U^2x-Ux}.\]
Hence the probability distribution does not converge unless $Ux=x$. We may also define the probability on a vertex to be the sum over the probabilities on its outgoing arcs. Using Kronecker's theorem, Aharonov et al \cite{Aharonov2000} showed that the probability  distribution over the vertices does not converge either. Nonetheless, the time-averaged probability distribution converges, under both notions.

In this section, we extend the discussion to the probability on any subset $S$ of the arcs. Suppose the arcs in $S$ are indexed by $\{\seq{j}{1}{2}{m}\}$. At time $k$, the entries of the following vector
\[\pmat{e_{j_1} & e_{j_2} & \cdots & e_{j_m}}^T U^kx\]
are the amplitudes on the arcs in $S$. Let
\[D_S = \pmat{e_{j_1} & e_{j_2} & \cdots & e_{j_m}}\pmat{e_{j_1} & e_{j_2} & \cdots & e_{j_m}}^T\]
be the characteristic matrix of $S$, that is, the diagonal matrix whose $(a,a)$-entry is $1$ if $a$ is in $S$. Then given initial state $x$, the probability that the quantum walk is on $S$ at time $k$ is
\[P_{x,S}(k) := x^* (U^k)^* D_S U^k x.\]
While this probability does not converge as $k$ goes to infinity, the \textsl{average probability} over time
\[\frac{1}{K} \sum_{k=0}^{K-1} P_{x,S}(k)\]
does as $K$ goes to infinity, for any initial state $x$. Before we prove it, we note that every unitary matrix $U$ can be written as
\[U = \sum_r e^{i\theta_r} F_r\]
where $e^{i\theta_r}$ is an eigenvalue of $U$, and $F_r$ is a Hermitian matrix representing the projection onto the eigenspace of $e^{i\theta_r}$.  This is usually called the \textsl{spectral decomposition} of $U$.

\begin{theorem}\label{lim_xS}
	Let $\seq{F}{1}{2}{m}$ be the spectral idempotents of the transition matrix $U$. Let $x$ be the initial state. For any subset $S$ of the arcs, the average probability that the quantum walk is on some arc of $S$ converges to
	\[\sum_r x^* F_r D_S F_r x.\]
	\label{limit}
\end{theorem}
\proof
Consider the spectral decomposition of $U$
\[U=\sum_r e^{i\theta_r} F_r.\]
It suffices to show that 
\[\frac{1}{K} \sum_{k=0}^{K-1} (U^k)^*D_SU^k\]
converges to 
\[\sum_r F_rD_S F_r\]
as $K$ goes to infinity. We have
\begin{align*}
(U^k)^* D_S U^k & =\left(\sum_r e^{-ik\theta_r}F_r\right)D_S \left(\sum_s e^{ik\theta_s}F_r\right)\\
&=\sum_r F_r D_S F_r+ \sum_{r\ne s} e^{ik(\theta_s-\theta_r)}F_rD_SF_s.
\end{align*}
Note that for all $r$ and $s$, the entries in $F_r D_S F_r$ and $F_rD_SF_s$ are constants, and remain unchanged when we take the average and the limit. Further
\[\frac{1}{K}\abs{\sum_{k=0}^{K-1} e^{ik(\theta_s-\theta_r)}} = \frac{1}{K} \abs{\frac{1 - e^{iK(\theta_s-\theta_r)}}{1-e^{i(\theta_s-\theta_r)}}} \le \frac{1}{K} \frac{2}{\abs{1-e^{i(\theta_s-\theta_r)}}}\]
which converges to zero as $K$ goes to infinity. Hence the only term that survives in
\[\lim_{K\to \infty}\frac{1}{K} \sum_{k=0}^{K-1} (U^k)^*D_SU^k\]
 is
\[\sum_r F_r D_S F_r. \tag*{\sqr53}\]

\section{Channels}
From the previous section, we see that the positive semidefinite matrices $F_r D_S F_r$ are crucial in computing the average limiting distribution of a quantum walk. Here is another way to interpret the relation between these matrices and the limit. The spectral idempotents $F_r$ satisfy that
\[\sum_r F_r^* F_r = I,\]
so the mapping on density matrices $\rho$ given by
\[\rho \mapsto \sum_r F_r \rho F_r^*\]
is a quantum channel. As $K$ goes to infinity, the average state
\[\frac{1}{K} \sum_{k=0}^{K-1} (U^k)^*\rho U^k\]
converges to the density matrix after the channel:
\[\sum_r F_r \rho F_r.\]
Now,  let $\abs{S}$ be the size of the subset $S$ and consider
\[\rho := \frac{1}{\abs{S}} D_S = \frac{1}{\abs{S}} e_{j_1}e_{j_1}^T + \cdots \frac{1}{\abs{S}} e_{j_m}e_{j_m}^T.\]
This is a density matrix of a mixed state---a convex combination of the pure states associated with the arcs in $S$. Hence, if the quantum walk starts with state $x$, then the average probability of the walk being on $S$ is the inner product of the density matrices $xx^*$ and $\sum_r F_r \rho F_r$, scaled by the size of $S$:
\[\abs{S}\Big<xx^*, \sum_r F_r  \rho F_r\Big> = \sum_r x^* F_r D_S F_r x.\]

\section{Mixing Times}
The mixing time of a quantum walk measures how far the average probability distribution is from the stationary distribution. Since the average probability on every subset converges, it does not make a big difference whether the distribution is on the arcs or on the vertices. In \cite{Aharonov2000}, Aharonov et al studied the mixing time on the vertices and obtained an upper bound for a general graph. They further showed that the mixing time of a quantum walk on an $n$-cycle with the Hadamard coin is bounded above by $O(n\log n)$, giving a quadratic speedup over the classical walk. We now extend some of their results to the arc probabilities.

Given $\epsilon$, define the \textsl{mixing time} $M_{\epsilon}$ to be the smallest $K$ such that for all $L\ge K$ and all initial states, the average distribution over $L$ steps on the arcs is $\epsilon$-close to the limiting distribution. More explicitly, let $P_j(k)$ denote the probability on the $j$-th arc, and let
\[D_j=e_je_j^T\]
be its characteristic matrix. The mixing time is is the minimum positive integer $K$ such that for all time $L\ge K$ and all initial states $x$, 
\[ \sum_j  \Bigg|\frac{1}{L}\sum_{k=0}^{L-1}P_{x,j}(k)- \sum_r x^* F_r D_j F_r x\Bigg| \le \epsilon.\]

\begin{lemma}\label{bound}
For a quantum walk with spectral decomposition
\[U = \sum_r \lambda_r F_r,\]
we have
\[\sum_j \Bigg|\frac{1}{K} \sum_{k=0}^{K-1} P_{x,j}(k)- \sum_r x^* F_r D_j F_r x\Bigg|\le \frac{2}{K} \sum_{r\ne s} \sum_j\frac{\sqrt{(F_r)_{jj} (F_s)_{jj}}}{\abs{\lambda_r - \lambda_s}}.\]
\end{lemma}

\proof
First note that for any $r$ and $s$, 
\begin{align*}
\abs{x^* F_r  D_j F_s x}
&= \abs{\grp{F_r e_j, x}\grp{F_s e_j, x}}\\
&\le \sqrt{(F_r)_{jj}} \norm{x} \sqrt{(F_s)_{jj}} \norm{x}\\
&= \sqrt{(F_r)_{jj} (F_s)_{jj}}.
\end{align*}
Let $\lambda_r = e^{i\theta_r}$ for some $\theta_r$. By Theorem \ref{limit}, for the $j$-th arc, 
\begin{align*}
\Bigg|\frac{1}{K} \sum_{k=0}^{K-1} P_{x,j}(k) - \sum_r x^* F_r D_j F_r x\Bigg|
&= \Bigg| \frac{1}{K}\sum_{k=0}^{K-1} \sum_{r\ne s} e^{ik(\theta_s-\theta_r)} x^* F_r D_j F_s x \Bigg|\\
&=\frac{1}{K} \Bigg| \sum_{r\ne s} \left(\sum_{k=0}^{K-1} e^{ik(\theta_s-\theta_r)}\right) x^*F_r D_j F_s x \Bigg|\\
&=\frac{1}{K} \Bigg| \sum_{r\ne s} \frac{1- e^{iK(\theta_s - \theta_r)}}{1-e^{i(\theta_s-\theta_r)}} x^*F_r D_j F_s x\Bigg|\\
&\le \frac{1}{K} \sum_{r\ne s}  \Bigg|  \frac{1- e^{iK(\theta_s - \theta_r)}}{1-e^{i(\theta_s-\theta_r)}}\Bigg| \abs{x^*F_r D_j F_s x}\\
&\le \frac{1}{K} \sum_{r\ne s} \frac{2}{\abs{e^{i\theta_s} - e^{i\theta_r}}} \abs{x^*F_r D_j F_s x}\\
&\le \frac{2}{K} \sum_{r \ne s} \frac{\sqrt{(F_r)_{jj} (F_s)_{jj}}}{\abs{\lambda_r - \lambda_s}}.
\end{align*}
Summing over all arcs yields the result.
\qed

One immediate consequence is that we can bound the mixing time of a quantum walk by its eigenvalue differences. This is an analogy to Lemma 4.3 in Aharonov \cite{Aharonov2000}.

\begin{corollary}
For a $\ell\times \ell$ transition matrix $U$ with spectral decomposition
\[U = \sum_r \lambda_r F_r,\]
we have 
\[M_{\epsilon} \le  \frac{2\ell}{\epsilon} \sum_{r\ne s} \frac{1}{\abs{\lambda_r - \lambda_s}}.\]
\end{corollary}

\proof
Since 
\[\sum_r F_r = I,\]
for all $r$ and $j$ we have
\[0\le (F_r)_{jj} \le 1.\]
Lemma \ref{bound} reduces to
\[\sum_j \Bigg|\frac{1}{K} \sum_{k=0}^{K-1} P_{x,j}(k)- \sum_r x^* F_r D_j F_r x\Bigg| \le \frac{2\ell}{K} \sum_{r\ne s} \frac{1}{\abs{\lambda_r - \lambda_s}}.\]
Thus for all $K$ such that 
\[K\ge \frac{2\ell}{\epsilon} \sum_{r\ne s} \frac{1}{\abs{\lambda_r - \lambda_s}},\]
the right hand side is bounded above by $\epsilon$.
\qed

\section{Average Mixing Matrix \label{avgmix}}
Let $U$ be the transition matrix of a discrete-time quantum walk. In this section, we consider the probabilities on the arcs only. They are given by the entries of 
\[U^k \circ \comp{U^k}.\]
Theorem \ref{limit} implies that the average probability distribution
\[\frac{1}{K} \sum_{k=0}^K U^k \circ \comp{U^k}\]
converges to the following matrix
\[\sum_r F_r \circ \comp{F_r}.\]
Note that this is a real symmetric matrix. Following Godsil's notion for continuous-time quantum walks \cite{Godsil2011}, we will denote this limit by $\widehat{M}$, and refer to it as the \textsl{average mixing matrix}. In \cite{Godsil2011}, Godsil established several properties of the continuous-time average mixing. We extend some of his results to discrete-time quantum walks.

The first observation is that $\widehat{M}$ is doubly-stochastic, in both the continuous and discrete cases. Thus each column of $\widehat{M}$ represents a probability distribution. For the $j$-th column of $\widehat{M}$, we define its \textsl{entropy} to be the negative expectation of the logarithm of its entries, that is,
\[-\sum_{\ell}\widehat{M}_{\ell j} \log\left(\widehat{M}_{\ell j}\right)\]
The entropy reaches its maximum if and only if the probability distribution with respect to the initial state $e_j$ is uniform. We also define the \textsl{total entropy} to be the sum of the entropies over all columns, although it is not an entropy. This invariant has been applied in structural pattern recognition. For example, Bai et al \cite{Bai2016} proposed a graph signature based on the entropy of the average mixing matrix of a graph. According to the experimental results, this entropic measure allows us to distinguish different structures.

In the continuous case, the average mixing matrix is positive semidefinite with eigenvalues no greater than one \cite{Godsil2011}.  We show that the same statement holds for the discrete average mixing matrix.

\begin{lemma}
The average mixing matrix $\widehat{M}$ is positive semidefinite, and its eigenvalues lie in $[0,1]$.
\end{lemma}

\proof
Since $F_r$ is positive semidefinite, its complex conjugate $F_r$ is also positive semidefinite. As a principal submatrix of $F_r \otimes \comp{F_r}$, the Schur product $F_r\circ \comp{F_s}$ is positive semidefinite. Hence the eigenvalues of $\widehat{M}$ are non-negative. It follows from
\[I= I\circ I=\left(\sum_r F_r\right) \left(\sum_s \comp{F_s}\right)=\widehat{M}+\sum_{r\ne s}F_r \circ \comp{F_s}\]
that the eigenvalues of $\widehat{M}$ are at most 1. On the other hand, since $\widehat{M}$ is doubly stochastic, $\one$ is an eigenvector of $\widehat{M}$ for the eigenvalue $1$.
\qed

If all the entries of $\widehat{M}$ are equal, we say the quantum walk has \textsl{uniform average mixing}. For example, the shunt-decomposition walk on $K_{3,3}$ with linear orders
\begin{gather*}f_1(1)=4,\quad f_1(2)=5,\quad f_1(3)=6;\\
	f_2(1)=5,\quad f_2(2)=6,\quad f_2(3)=4;\\
	f_3(1)=6,\quad f_3(2)=4,\quad f_3(3)=5;\\
	f_4(1)=2,\quad f_4(2)=1,\quad f_4(3)=3;\\
	f_5(1)=3,\quad f_5(2)=2,\quad f_5(3)=1;\\
	f_6(1)=1,\quad f_6(2)=3,\quad f_6(3)=2,
\end{gather*}
and the following coin
\[C = \frac{1}{7}\pmat{
	-2 & 3 & 6\\ 
	6 & -2 & 3\\
	3 & 6 & -2}\]
admits uniform average mixing. It indicates that in the limit, regardless of the initial arc, the quantum walk may be found on any arc with the same probability. Our next goal is to establish necessary and sufficient conditions for uniform average mixing to occur.

Suppose the spectral decomposition of $U$ is
\[U = \sum_r \theta_r F_r.\]
We say $U$ is \textsl{walk-regular} if $F_r$ has constant diagonal for each $r$. This is related, but not equivalent, to the walk-regularity of a graph; the latter plays a role in studying the average mixing matrix of a continuous-time quantum walk.

\begin{lemma}
Let $U$ be an $\ell\times \ell$ unitary matrix. If $m_r$ is the multiplicity of the $r$-th eigenvalue of $U$, then 
\[\tr(\widehat{M})\ge \frac{1}{\ell} \sum_r m_r^2.\]
Further, equality holds if and only if $U$ is walk-regular.
\label{tr}
\end{lemma}
\proof
Consider the spectral decomposition 
\[U = \sum_r \lambda_r F_r.\]
Since $F_r$ is positive semidefinite, its diagonal entries are non-negative and
\[\tr(F_r) = m_r.\]
By the Cauchy-Schwarz inequality,
\[\tr(F_r \circ \comp{F_r}) \ge  \frac{1}{\ell} \tr(F_r)^2 = \frac{1}{\ell} m_r^2.\]
Hence
\[\tr(\widehat{M}) \ge \frac{1}{\ell} \sum_r m_r ^2. \]
Equality holds if and only if each $F_r$ has constant diagonal $m_r/\ell $.
\qed

\begin{corollary}\label{secd}
For any unitary matrix $U$ we have $\tr(\widehat{M}) \ge 1$. Equality holds if and only if  $U$ is walk-regular and has simple eigenvalues.
\end{corollary}
\proof
The inequality follows from Lemma \ref{tr} and 
\[\sum_r m_r = \ell.\]
If equality holds, then $m_r=1$ for all $r$, and $U$ is walk-regular. 
\qed

\begin{theorem}
The following statements are equivalent.
\begin{enumerate}[(i)]
	\item The quantum walk admits uniform average mixing.
	\item $\tr(\widehat{M})=1$.
	\item $U$ is walk-regular with simple eigenvalues.
\end{enumerate}
\end{theorem}
\proof
Let $U$ be an $\ell \times \ell$ unitary matrix. To see that (i) implies (ii), note that if uniform average mixing occurs, then all entries of $\widehat{M}$ are $1/\ell$, so $\tr(\widehat{M})=1$. Corollary \ref{secd} shows that (ii) implies (iii). Now suppose (iii) holds. Then the spectral decomposition of $U$ is
\[U = \sum_r z_r z_r^*,\]
where $z_r$ is an normalized eigenvector of $U$ for the eigenvalue $\theta_r$. Since $U$ is walk-regular, for each $r$, the entries of $x_r$ have the same absolute value. Thus 
\[\widehat{M} = \sum_r (z_r z_r^*) \circ (z_r z_r^*) = \sum_r (z_r \circ z_r) (z_r \circ z_r)^*\]
is flat. Therefore (iii) implies (i).
\qed

The average mixing matrix $\widehat{M}$ records the limiting probability on each arc, given that the walk started with an arc. One may also compute 
\[\lim_{K\to \infty} \frac{1}{K} \sum_{k=0}^{K-1} P_{x, S}(k),\]
that is, the limiting probability on a set $S$ of arcs, given initial state $x$. In \cite{Aharonov2000}, Aharonov et al discussed when the limiting probability distribution is uniform over all vertices. We show that this is guaranteed whenever $\widehat{M}$ is flat, regardless of the initial state. 

\begin{lemma}
	Suppose $U$ has simple eigenvalues with corresponding eigenvectors $\seq{z}{1}{2}{\ell}$. Let 
	\[x = \sum_r a_r z_r,\]
	where $\sum_r \abs{a_r}^2=1$. Then the limiting probability that the quantum walk is on $S$, given initial state $x$, is
	\[\lim_{K\to \infty} \frac{1}{K} \sum_{k=0}^{K-1} P_{x, S}(k) = \sum_r \abs{a_r}^2 z_r^* D_S z_r.\]
	Moreover, if $U$ is walk-regular, then the limiting probability distribution is constant over $V(X)$, regardless of the initial state.
\end{lemma}
\proof
The first statement follows from Theorem \ref{lim_xS}. Now, if $U$ is walk-regular with simple eigenvalues, then each eigenvector $z_r$ is flat, and $z_r^* D_S z_r$ depends only on the size of $S$. 
\qed

Finally, we prove some algebraic properties of the average mixing matrix. They rely on the well-known fact that a commutative semisimple matrix algebra with identity has a basis of orthogonal idempotents. The following is a useful lemma in proving the existence of some phenomenon of quantum walks.

\begin{lemma}
If the entries of $U$ are algebraic over $\rats$, then the entries of $\widehat{M}$ are algebraic over $\rats$.
\end{lemma}
\proof
Suppose $U$ has algebraic entries. Then its eigenvalues are all algebraic. Let $\fld$ be the smallest field containing the eigenvalues of $U$. Let $\cB$ be the matrix algebra generated by $U$ over $\fld$. To show that $\cB$ is semisimple, pick $N\in \cB$ with $N^2=0$. Since $U$ is unitary, the algebra $\cB$ is closed under conjugate transpose and contains the identity. It follows from $(N^*)^2=0$ that
\begin{align*}
0&= \tr((N^*)^2N^2)\\
&= \tr(N^*N N^*N)\\
&= \tr((N^*N)^*(N^*N)).
\end{align*}
Thus $N^*N=0$. Applying the trace again to $N^*N$, we see that $N=0$. Therefore, the spectral idempotents $F_r$ of $U$ are polynomials in $U$ with algebraic coefficients. Hence the entries in
\[\widehat{M} = \sum_r F_r \circ \comp{F_r}\]
are algebraic over $\rats$.
\qed

In continuous-time quantum walks, rationality has been used in different ways to determine which graphs admit instantaneous uniform mixing \cite{Godsil2014a}. It is also known that the entries of a continuous average mixing matrix are all rational \cite{Godsil2011}. We show that the discrete average mixing matrix enjoys the same property, given that all entries of the transition matrix are rational.

\begin{lemma}
If the entries of $U$ are rational, then the entries of  $\widehat{M}$ are rational.
\end{lemma}

\proof
Let $\fld$ be the smallest field containing the eigenvalues of $U$. Let $\sigma$ be an automorphism of $\fld$. Since $U$ is rational, we have
\[U=U^{\sigma}=\sum_r \lambda_r^{\sigma} F_r^{\sigma}.\]
Moreover, since $\lambda_r^{\sigma}$ is also an eigenvalue of $U$, the set of idempotents $\{F_r\}$ is closed under field automorphisms. Thus
\[\widehat{M}=\sum_r F_r\circ F_r^T\]
is fixed by all automorphisms of $\fld$ and must be rational.
\qed

\section{Traces and Total Entropies}
In this section, we present some numerical evidence for how the structure of a graph affects the limiting distribution of a quantum walk. The graphs we test are all cubic on up to 8 vertices. One of the parameters of the walk we examine is $\tr(\widehat{M})$, as every diagonal entry of the average mixing matrix represents the probability of a walk returning to the arc it started with. Thus, a higher trace indicates a higher tendency for the walk to stay at home. For the arc-reversal model, we also compute the total entropies of $\widehat{M}$, which measure how close the limiting probability distribution is to uniform.

\subsection{Rotation System\label{data_arcrev}}
For each graph, we enumerate all the rotation systems, and investigate the associated arc-reversal walks. Since these graphs are all cubic, the following circulant with simple eigenvalues is a candidate for the coin:
\[C = \frac{1}{7}\pmat{
-2 & 3 & 6\\ 
6 & -2 & 3\\
3 & 6 & -2}.\]
We then compute the average mixing matrix and its trace for each rotation system. 

There are $2^n$ rotation systems for a cubic graph on $n$ vertices. The following table provides the number of rotation systems with the same genus and the same trace. For cubic graphs that do not have a name, we put their graph6 strings. It appears that, for the same graph, embeddings on higher-genus surfaces give lower traces.

\begin{longtable}{c||c||c||c}
\hline
\hline
graph & genus  & trace & number of rotation systems\\
\hline
\hline
\multirow{3}{*}{$K_4$}
& $0$ & $3.000000$ & $2$\\
\cline{2-4}
& $1$ & $1.753644$ & $8$\\
\cline{2-4}
& $1$ & $1.694295$ & $6$\\
\hline
\hline
\multirow{3}{*}{$K_{3,3}$}
& $1$ & $2.201010$ & $36$\\
\cline{2-4}
& $1$ & $2.111111$ & $4$\\
\cline{2-4}
& $2$ & $1.052644$ & $24$\\
\hline
\hline
\multirow{8}{*}{$K_2\square K_3$}
& $0$ & $3.255278$ & $2$\\
\cline{2-4}
& $1$ & $2.105870$ & $12$\\
\cline{2-4}
& $1$ & $2.089084$ & $6$\\
\cline{2-4}
& $1$ & $1.932964$ & $2$\\
\cline{2-4}
& $1$ & $1.918699$ & $12$ \\
\cline{2-4}
& $1$ & $1.866536$ & $6$\\
\cline{2-4}
& $2$ & $1.340085$ & $12$\\
\cline{2-4}
& $2$ & $1.187163$ & $12$\\
\hline
\hline
\end{longtable}
\newpage
\begin{longtable}{c||c||c||c}
\hline
\hline
graph & genus  & trace & number of rotation systems\\
\hline
\hline
\multirow{13}{*}{$Q_3$}
& $0$ & $4.500000$ & $2$\\
\cline{2-4}
& $1$ & $3.000000$ & $8$\\
\cline{2-4}
& $1$ & $2.744344$ & $16$\\
\cline{2-4}
& $1$ & $2.625302$ & $24$\\
\cline{2-4}
& $1$ & $2.446501$ & $6$\\
\cline{2-4}
& $2$ & $1.980844$ & $8$\\
\cline{2-4}
& $2$ & $1.746199$ & $48$\\
\cline{2-4}
& $2$ & $1.694728$ & $24$\\
\cline{2-4}
& $2$ & $1.694295$ & $8$\\
\cline{2-4}
& $2$ & $1.688522$ & $16$\\
\cline{2-4}
& $2$ & $1.680899$ & $48$\\
\cline{2-4}
& $2$ & $1.679098$ & $24$\\
\cline{2-4}
& $2$ & $1.524169$ & $24$\\
\hline
\hline
\multirow{17}{*}{$X(\ints_8,\{3,4,5\}$)}
 & $1$ & $2.520693$ & $16$\\\cline{2-4}
& $1$ & $2.513717$ & $8$\\\cline{2-4}
& $1$ & $2.447996$ & $8$\\\cline{2-4}
& $1$ & $2.077332$ & $16$\\\cline{2-4}
& $1$ & $2.00942$ & $8$\\\cline{2-4}
& $2$ & $1.914302$ & $6$\\\cline{2-4}
& $2$ & $1.886742$ & $16$\\\cline{2-4}
& $2$ & $1.752042$ & $32$\\\cline{2-4}
& $2$ & $1.658131$ & $2$\\\cline{2-4}
& $2$ & $1.650867$ & $16$\\\cline{2-4}
& $2$ & $1.615907$ & $16$\\\cline{2-4}
& $2$ & $1.599010$ & $16$\\\cline{2-4}
& $2$ & $1.598366$ & $16$\\\cline{2-4}
& $2$ & $1.586057$ & $16$\\\cline{2-4}
& $2$ & $1.566265$ & $16$\\\cline{2-4}
& $2$ & $1.557058$ & $32$\\\cline{2-4}
& $2$ & $1.460324$ & $16$\\
\hline
\hline
\end{longtable}
\newpage
\begin{longtable}{c||c||c||c}
\hline
\hline
graph & genus  & trace & number of rotation systems\\
\hline
\hline
\multirow{44}{*}{\texttt{GCZJd\_}}
& $0$ & $3.019811$ & $2$\\\cline{2-4}
& $1$ & $2.542603$ & $2$\\\cline{2-4}
& $1$ & $2.470562$ & $4$\\\cline{2-4}
& $1$ & $2.429130$ & $4$\\\cline{2-4}
& $1$ & $2.362862$ & $4$\\\cline{2-4}
& $1$ & $2.354247$ & $8$\\\cline{2-4}
& $1$ & $2.326463$ & $2$\\\cline{2-4}
& $1$ & $2.321465$ & $8$\\\cline{2-4}
& $1$ & $2.299543$ & $4$\\\cline{2-4}
& $1$ & $2.271257$ & $4$\\\cline{2-4}
& $1$ & $2.259590$ & $8$\\\cline{2-4}
& $1$ & $2.250146$ & $4$\\\cline{2-4}
& $1$ & $2.214541$ & $4$\\\cline{2-4}
& $1$ & $2.178556$ & $4$\\\cline{2-4}
& $1$ & $2.115562$ & $2$\\\cline{2-4}
& $1$ & $2.079810$ & $8$\\\cline{2-4}
& $1$ & $1.75710$ & $8$\\\cline{2-4}
& $2$ & $1.714738$ & $8$\\\cline{2-4}
& $2$ & $1.713346$ & $4$\\\cline{2-4}
& $2$ & $1.709795$ & $8$\\\cline{2-4}
& $2$ & $1.690956$ & $4$\\\cline{2-4}
& $2$ & $1.685911$ & $8$\\\cline{2-4}
& $2$ & $1.649712$ & $8$\\\cline{2-4}
& $2$ & $1.622273$ & $8$\\\cline{2-4}
& $2$ & $1.618072$ & $8$\\\cline{2-4}
& $2$ & $1.605598$ & $8$\\\cline{2-4}
& $2$ & $1.600225$ & $8$\\\cline{2-4}
& $2$ & $1.599082$ & $8$\\\cline{2-4}
& $2$ & $1.571617$ & $8$\\\cline{2-4}
& $2$ & $1.571597$ & $8$\\\cline{2-4}
& $2$ & $1.563497$ & $8$\\\cline{2-4}
 & $2$ & $1.560172$ & $8$\\\cline{2-4}
& $2$ & $1.553949$ & $4$\\\cline{2-4}
 & $2$ & $1.541723$ & $4$\\\cline{2-4}
& $2$ & $1.535939$ & $8$\\\cline{2-4}
 & $2$ & $1.531083$ & $8$\\\cline{2-4}
& $2$ & $1.521709$ & $4$\\\cline{2-4}
& $2$ & $1.505587$ & $4$\\\cline{2-4}
& $2$ & $1.505097$ & $8$\\\cline{2-4}
& $2$ & $1.486416$ & $4$\\\cline{2-4}
& $2$ & $1.468886$ & $4$\\\cline{2-4}
& $2$ & $1.466287$ & $4$\\\cline{2-4}
& $2$ & $1.457042$ & $8$\\\cline{2-4}
& $2$ & $1.421791$ & $4$\\
\hline
\hline
\multirow{19}{*}{\texttt{GCXmd\_} }
 & $0$ & $3.681253$ & $4$\\\cline{2-4}
& $1$ & $2.560923$ & $16$\\\cline{2-4}
& $1$ & $2.553190$ & $4$\\\cline{2-4}
& $1$ & $2.535164$ & $16$\\\cline{2-4}
& $1$ & $2.502552$ & $8$\\\cline{2-4}
& $1$ & $2.351188$ & $16$\\\cline{2-4}
& $1$ & $2.331725$ & $16$\\\cline{2-4}
& $1$ & $2.278575$ & $8$\\\cline{2-4}
& $1$ & $2.150255$ & $8$\\\cline{2-4}
& $2$ & $1.769929$ & $16$\\\cline{2-4}
& $2$ & $1.746908$ & $16$\\\cline{2-4}
& $2$ & $1.742805$ & $32$\\\cline{2-4}
& $2$ & $1.660490$ & $32$\\\cline{2-4}
& $2$ & $1.619199$ & $8$\\\cline{2-4}
& $2$ & $1.544993$ & $16$\\\cline{2-4}
& $2$ & $1.541696$ & $8$\\\cline{2-4}
& $2$ & $1.534019$ & $16$\\\cline{2-4}
& $2$ & $1.466143$ & $8$\\\cline{2-4}
& $2$ & $1.450838$ & $8$\\
\hline
\hline
\end{longtable}
\newpage
\begin{longtable}{c||c||c||c}
\hline
\hline
graph & genus  & trace & number of rotation systems\\
\hline
\hline
\multirow{16}{*}{\texttt{GCY\^{}B\_}}
& $1$ & $2.402758$ & $24$\\\cline{2-4}
& $1$ & $2.351016$ & $24$\\\cline{2-4}
& $1$ & $2.263351$ & $12$\\\cline{2-4}
& $1$ & $2.018540$ & $4$\\\cline{2-4}
& $2$ & $1.849363$ & $12$\\\cline{2-4}
& $2$ & $1.805881$ & $24$\\\cline{2-4}
& $2$ & $1.775098$ & $24$\\\cline{2-4}
& $2$ & $1.766025$ & $24$\\\cline{2-4}
& $2$ & $1.736332$ & $12$\\\cline{2-4}
& $2$ & $1.714296$ & $8$\\\cline{2-4}
& $2$ & $1.699936$ & $4$\\\cline{2-4}
& $2$ & $1.676129$ & $24$\\\cline{2-4}
& $2$ & $1.656324$ & $12$\\\cline{2-4}
& $2$ & $1.562093$ & $12$\\\cline{2-4}
& $2$ & $1.519301$ & $12$\\\cline{2-4}
& $2$ & $1.455918$ & $24$\\
\hline
\hline
\caption{$\tr(\widehat{M})$ of arc-reversal walks based on rotation systems}
\end{longtable}

Although there are $2^n$ rotation systems in total, many of them are equivalent, in the sense that there is an isomorphism between the associated gem that preserves each edge coloring. As indicated by the following table, the total entropy seems to distinguish non-equivalent embeddings, and for most of the time it increases as the genus increases.

\begin{longtable}{c||c||c||c}
\hline
\hline
graph & rotation system & genus & entropy\\
\hline
\hline
\multirow{3}{*}{\includegraphics[width=1.3cm]{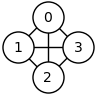}}
& $0: (1, 2, 3), 1: (0, 3, 2), 2: (0, 1, 3), 3: (0, 2, 1)$ & $0$ & $25.364055$\\
\cline{2-4}
& $0: (1, 2, 3), 1: (0, 2, 3), 2: (0, 1, 3), 3: (0, 1, 2)$ & $1$ & $27.490608$\\
\cline{2-4}
& $0: (1, 2, 3), 1: (0, 3, 2), 2: (0, 1, 3), 3: (0, 1, 2)$ & $1$ & $27.763049$ \\
\hline\hline
\multirow{5}{*}{\includegraphics[width=2.8cm]{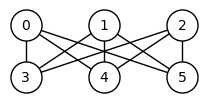}}
&\pbox{20cm}{$0:(3,5,4),1:(3,5,4),2:(3,4,5),$\\$3:(0,2,1),4:(0,1,2),5:(0,2,1)$} & $1$ & $47.42653$\\
 \cline{2-4} 
& \pbox{20cm}{$0:(3,5,4),1:(3,5,4),2:(3,5,4),$\\$3:(0,2,1),4:(0,2,1),5:(0,2,1)$} & $1$ & $47.862470$\\
 \cline{2-4} 
& \pbox{20cm}{$0: (4, 3, 5), 1: (4, 3, 5), 2: (4, 3, 5),$\\$ 3: (1, 0, 2), 4: (0, 1, 2), 5: (0, 1, 2)$} & $2$ & $52.001943$\\
 \hline \hline
 \end{longtable}
 \newpage
\begin{longtable}{c||c||c||c}
\hline
\hline
graph & rotation system & genus & entropy\\
\hline
\hline
\multirow{13}{*}{\includegraphics[width=2.8cm]{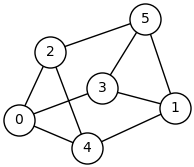}} 
& \pbox{20cm}{$0: (2, 3, 4),
  1: (4, 3, 5),
  2: (0, 4, 5),$\\$
  3: (1, 0, 5),
  4: (0, 1, 2),
  5: (2, 1, 3)$} & $0$ & $45.68992$\\
 \cline{2-4}  & \pbox{20cm}{$ 0: (2, 3, 4),
        1: (3, 4, 5),
        2: (4, 0, 5),$\\$
        3: (0, 1, 5),
        4: (0, 1, 2),
        5: (2, 1, 3)$} & $1$ & $48.861877$\\
 \cline{2-4} & \pbox{20cm}{$0: (2, 3, 4),
    1: (3, 4, 5),
    2: (0, 4, 5),$\\$
    3: (0, 1, 5),
    4: (0, 1, 2),
    5: (1, 2, 3)$} & $1$ & $48.864165$\\
 \cline{2-4} & \pbox{20cm}{$0: (2, 3, 4),
      1: (3, 4, 5),
      2: (0, 4, 5),$\\$
      3: (1, 0, 5),
      4: (0, 1, 2),
      5: (2, 1, 3)$} & $1$ & $48.981188$\\
 \cline{2-4} & \pbox{20cm}{$0: (2, 3, 4),
     1: (3, 4, 5),
     2: (0, 4, 5),$\\$
     3: (1, 0, 5),
     4: (1, 0, 2),
     5: (2, 1, 3)$} & $1$ & $49.081280$\\
 \cline{2-4} & \pbox{20cm}{$0: (2, 3, 4),
   1: (3, 4, 5),
   2: (0, 4, 5),$\\$
   3: (0, 1, 5),
   4: (0, 1, 2),
   5: (2, 1, 3)$} & $1$ & $49.692402$\\
 \cline{2-4} & \pbox{20cm}{$0: (2, 3, 4),
  1: (3, 4, 5),
  2: (0, 4, 5),$\\$
  3: (0, 1, 5),
  4: (1, 0, 2),
  5: (2, 1, 3)$} & $2$ & $51.692584$\\
 \cline{2-4}  & \pbox{20cm}{$0: (2, 3, 4),
  1: (3, 4, 5),
  2: (4, 0, 5),$\\$
  3: (0, 1, 5),
  4: (1, 0, 2),
  5: (2, 1, 3)$} & $2$ & $51.928224$\\
  \hline \hline 
  \multirow{20}{*}{\includegraphics[width=2.8cm]{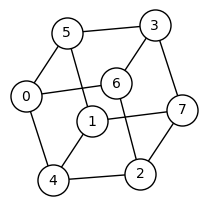}}  
& \pbox{20cm}{$0: (4, 5, 6),
   1: (5, 4, 7),
   2: (4, 6, 7),
   3: (6, 5, 7),$\\$
   4: (1, 0, 2),
   5: (0, 1, 3),
   6: (2, 0, 3),
   7: (1, 2, 3)$} & 0 & $62.411249$\\
 \cline{2-4}   & \pbox{20cm}{$0: (4, 5, 6),
    1: (5, 4, 7),
    2: (4, 6, 7),
    3: (6, 5, 7),$\\$
    4: (0, 1, 2),
    5: (1, 0, 3),
    6: (0, 2, 3),
    7: (2, 1, 3)$} & $1$ & $67.363643$\\
 \cline{2-4}  & \pbox{20cm}{$0: (4, 5, 6),
     1: (4, 5, 7),
     2: (4, 6, 7),
     3: (5, 6, 7),$\\$
     4: (0, 1, 2),
     5: (0, 1, 3),
     6: (0, 2, 3),
     7: (1, 2, 3)$} & $1$ & $68.359584$\\
 \cline{2-4}  & \pbox{20cm}{$0: (4, 5, 6),
      1: (4, 5, 7),
      2: (4, 6, 7),
      3: (6, 5, 7),$\\$
      4: (1, 0, 2),
      5: (0, 1, 3),
      6: (2, 0, 3),
      7: (1, 2, 3)$} & $1$ & $69.625653$\\
 \cline{2-4}  & \pbox{20cm}{$0: (4, 5, 6),
       1: (4, 5, 7),
       2: (4, 6, 7),
       3: (6, 5, 7),$\\$
       4: (0, 1, 2),
       5: (0, 1, 3),
       6: (2, 0, 3),
       7: (1, 2, 3)$} & $1$ & $69.919303$\\
 \cline{2-4}   & \pbox{20cm}{$0: (4, 5, 6),
        1: (4, 5, 7),
        2: (4, 6, 7),
        3: (5, 6, 7),$\\$
        4: (1, 0, 2),
        5: (1, 0, 3),
        6: (2, 0, 3),
        7: (2, 1, 3)$} & $1$ & $70.331970$\\
 \cline{2-4}  & \pbox{20cm}{$0: (4, 5, 6),
         1: (4, 5, 7),
         2: (6, 4, 7),
         3: (5, 6, 7),$\\$
         4: (0, 1, 2),
         5: (1, 0, 3),
         6: (0, 2, 3),
         7: (1, 2, 3)$} & $2$ & $71.633476$\\
 \cline{2-4}  & \pbox{20cm}{$0: (4, 5, 6),
        1: (5, 4, 7),
        2: (6, 4, 7),
        3: (6, 5, 7),$\\$
        4: (0, 1, 2),
        5: (0, 1, 3),
        6: (0, 2, 3),
        7: (2, 1, 3)$} & $2$ & $72.161631$\\
 \cline{2-4}  & \pbox{20cm}{$0: (4, 5, 6),
         1: (4, 5, 7),
         2: (6, 4, 7),
         3: (5, 6, 7),$\\$
         4: (0, 1, 2),
         5: (0, 1, 3),
         6: (2, 0, 3),
         7: (1, 2, 3)$} &  $2$ & $72.270262$ \\
 \cline{2-4}   & \pbox{20cm}{$0: (4, 5, 6),
          1: (5, 4, 7),
          2: (6, 4, 7),
          3: (6, 5, 7),$\\$
          4: (0, 1, 2),
          5: (1, 0, 3),
          6: (0, 2, 3),
          7: (2, 1, 3)$} & $2$ & $72.758589$\\
 \cline{2-4}  & \pbox{20cm}{$0: (4, 5, 6),
           1: (4, 5, 7),
           2: (4, 6, 7),
           3: (5, 6, 7),$\\$
           4: (0, 1, 2),
           5: (1, 0, 3),
           6: (0, 2, 3),
           7: (2, 1, 3)$} & $2$ & $73.548157$\\
 \cline{2-4}  & \pbox{20cm}{$0: (4, 5, 6),
            1: (4, 5, 7),
            2: (4, 6, 7),
            3: (5, 6, 7),$\\$
            4: (0, 1, 2),
            5: (1, 0, 3),
            6: (2, 0, 3),
            7: (2, 1, 3)$} & $2$ & $73.595086$ \\
 \cline{2-4}  & \pbox{20cm}{$0: (4, 5, 6),
             1: (5, 4, 7),
             2: (6, 4, 7),
             3: (5, 6, 7),$\\$
             4: (0, 1, 2),
             5: (0, 1, 3),
             6: (0, 2, 3),
             7: (1, 2, 3)$} & $2$ & $75.176378$ \label{entropy}    \\
 \hline \hline       
 \caption{Total entropy of $\widehat{M}$ of arc-reversal walks based on rotation systems}
\end{longtable}

\subsection{Shunt-Decompositions \label{data_shunt}}
We now switch to the shunt-decomposition model with the following $d\times d$ coin:
\[B = \left(\frac{1}{\sqrt{d}} e^{2(j-k)^2\pi i /d}\right)_{jk}.\]
This is a unitary circulant for odd $d$. When $d=3$, it is self-congruent under all permutations, so we do not need to specify an ordering of the shunts to run a quantum walk. Again, for each graph, we enumerate all shunt-decompositions, and compute the average mixing matrices and their traces. 

The following table lists only one shunt-decomposition for each cycle structure. The data indicates that for the same graph, the symmetric shunt-decompositions, if any, always gives the highest trace.

\begin{longtable}{c||c||c}
\hline \hline
graph & shunt-decomposition  & trace\\
\hline \hline
\multirow{2}{*}{\includegraphics[width=1cm]{K4}}& $\{(0,1)(2,3), (0,2,1,3), (0,3,1,2)\}$ &  $2.3333$\\\cline{2-3}
& $\{(0,1)(2,3), (0,2)(1,3), (0,3)(1,2)\}$ &  $2.6667$\\
\hline\hline
\multirow{4}{*}{\includegraphics[width=2.5cm]{K33}}
& $\{(0,3)(1,5,2,4), (0,4,2,5)(1,3), (0,5,1,4)(2,3)\}$ & $2.3333$\\
\cline{2-3}
& $\{(0,3)(1,4)(2,5), (0,4,2,3,1,5), (0,5,1,3,2,4)\}$ & $2.3333$\\
\cline{2-3}
& $\{(0,3,2,5,1,4), (0,4,2,3,1,5), (0,5,2,4,1,3)\}$ & $3.6667$\\
\cline{2-3}
& $\{(0,3)(1,5)(2,4), (0,4)(1,3)(2,5), (0,5)(1,4)(2,3)\}$ & $3.6667$ \\
\hline \hline 
\multirow{4}{*}{\includegraphics[width=2.1cm]{K2boxK3}}
& $\{(0,2,4)(1,3,5), (0,3)(1,4,2,5), (0,4,1,3,5,2)\}$ & $1.6630$ \\
\cline{2-3}
& $\{(0,2)(1,4)(3,5), (0,3,1,5,2,4), (0,4,2,5,1,3)\}$ & $1.7482$\\
\cline{2-3}
& $\{(0,2,4)(1,3,5), (0,4,2)(1,5,3), (0,3)(1,4)(2,5)\}$ & $2.3665$\\
\cline{2-3}
& $\{(0,2)(1,4)(3,5), (0,3)(1,5)(2,4), (0,4)(1,3)(2,5)\}$ & $2.6458$\\
\hline \hline
\multirow{6}{*}{\includegraphics[width=2.5cm]{Q3}}
& $\{(0,4,1,5)(2,7)(3,6), (0,5,3,7,1,4)(2,6), (0,6)(2,4)(1,5,3,7)\}$ & $1.7987$\\
\cline{2-3}
& $\{(0,4,2,6)(1,5,3,7), (0,5,1,4)(2,7)(3,6), (0,6,2,4,1,7,3,5)\}$ & $1.8738$\\
\cline{2-3}
& $\{(0,4,1,5)(2,6,3,7), (0,5,1,7,3,6)(2,4), (0,6,2,7,1,4)(3,5)\}$ & $1.9065$\\
\cline{2-3}
& $\{(0,4,1,7,2,6,3,5), (0,5,3,6,2,7,1,4), (0,6)(1,5)(2,4)(3,7)\}$ & $2.3241$\\
\cline{2-3}
& $\{(0,4,2,6)(1,7,3,5), (0,5,3,6,2,7,1,4), (0,6,3,7,2,4,1,5)\}$ & $2.4939$ \\
\cline{2-3}
& $\{(0,4)(1,5)(2,6)(3,7), (0,5,3,6)(1,7)(2,4), (0,6,3,5)(1,4)(2,7)\}$ & $2.5232$ \\
\cline{2-3}
& $\{(0,4)(1,5)(2,6)(3,7), (0,5,3,6)(1,4,2,7), (0,6,3,5)(1,7,2,4)\}$ & $2.6670$\\
\cline{2-3}
& $\{(0,4,1,5)(2,6,3,7), (0,5,3,6)(1,4,2,7), (0,6,2,4)(1,7,3,5)\}$ & $3.1667$ \\
\cline{2-3}
& $\{(0,4)(1,5)(2,7)(3,6), (0,5)(1,4)(2,6)(3,7), (0,6)(1,7)(2,4)(3,5)\}$ & $3.6667$\\
\hline \hline
\multirow{9}{*}{\includegraphics[width=2.5cm]{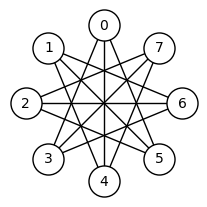}}
& $\{(0,3)(1,5,2,6)(4,7), (0,4)(1,6,3,7,2,5)\}$ & $1.8266$\\
\cline{2-3}
& $\{(0,3,7,4,1,6,2,5),(0,4)(1,5)(2,7)(3,6),(0,5,2,6,1,4,7,3)\}$ & $1.829313$\\
\cline{2-3}
 & $\{(0,3,7,4)(1,5,2,6),(0,4,7,2,5,1,6,3),(0,5)(1,4)(2,7,3,6)\}$ & $1.888459$\\
\cline{2-3}
& $\{(0,3)(1,6)(2,5)(4,7), (0,4,1,5)(2,6)(3,7), (0,3)(1,6)(2,5)(4,7)\}$ & $2.098040$\\
\cline{2-3}
& $\{(0,3,6,1,4,7,2,5),(0,4)(1,5)(2,6)(3,7),(0,5,2,7,4,1,6,3)\}$ & $2.218323$\\
\cline{2-3}
& $\{(0,3,6,2,7,4,1,5),(0,4)(1,6,3,7,2,5), (0,5,2,6,1,4,7,3)\}$ & $2.286421$\\
\cline{2-3}
& $\{(0,3)(1,6)(2,5)(4,7),(0,4,1,5)(2,6,3,7),(0,5,1,4)(2,7,3,6)\}$ & $2.4376$\\
\cline{2-3}
& $\{(0,3)(1,6)(2,5)(4,7), (0,4)(1,5)(2,6)(3,7), (0,5)(1,4)(2,7)(3,6)\}$ & $3.045749$\\
\hline \hline
\multirow{9}{*}{\includegraphics[width=2.5cm]{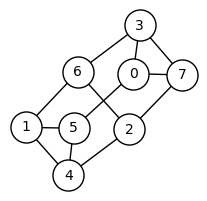}}
& $\{(0,3)(2,6,1,5)(4,7), (0,4)(1,6,3,7,2,5), (0,5)(1,4)(2,7,3,6)\}$ & $1.826585$\\
\cline{2-3}
& $\{((0,3,7,4,1,6,2,5),(0,4)(1,5)(2,7)(3,6),(0,5,2,6,1,4,7,3)\}$ & $1.829313$\\
\cline{2-3}
& $\{(0,3,7,4)(1,5,2,6),(0,4,7,2,5,1,6,3),(0,5)(1,4)(2,7,3,6)\}$ & $1.888459$\\
\cline{2-3}
& $\{(0,3)(1,6)(2,5)(4,7), (0,4,1,5)(2,6)(3,7), (0,5,1,4)(2,7)(3,6)\}$ & $2.098040$\\
\cline{2-3}
& $\{(0,3,6,1,4,7,2,5), (0,4)(1,5)(2,6)(3,7), (0,5,2,7,4,1,6,3)\}$ & $2.218323$\\
\cline{2-3}
& $\{(0,3,6,2,7,4,1,5),(0,4)(1,6,3,7,2,5)\}$ & $2.286421$\\
\cline{2-3}
& $\{(0,3)(1,6)(2,5)(4,7), (0,4,1,5)(2,6,3,7),(0,5,1,4)(2,7,3,6)\}$& $2.437614$\\
\cline{2-3}
& $\{(0,3)(1,6)(2,5)(4,7), (0,4)(1,5)(2,6)(3,7), (0,5)(1,4)(2,7)(3,6)\}$ & $3.045749$\\
\hline\hline
\multirow{5}{*}{\includegraphics[width=2.5cm]{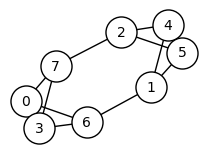}}
& $\{(0,3,6)(2,7)(1,4,5), (0,6,3,7)(1,5,2,4), (0,7,3)(1,6)(2,5,4)\}$ & $1.769242$\\
\cline{2-3}
& $\{(0,3,7,2,4,5,1,6), (0,6,1,5,4,2,7,3), (0,7)(1,4)(2,5)(3,6)\}$ & $2.186914$ \\
\cline{2-3}
& $\{(0,3)(1,6)(2,7)(4,5), (0,6,3,7)(1,4)(2,5), (0,7,3,6)(1,5)(2,4)\}$ & $2.676998$\\
\cline{2-3}
& $\{(0,3)(1,6)(2,7)(4,5), (0,6,3,7)(1,4,2,5), (0,7,3,6)(1,5,2,4)\}$ & $2.983961$\\
\cline{2-3}
& $\{(0,3)(1,6)(2,7)(4,5), (0,6)(1,4)(2,5)(3,7), (0,7)(1,5)(2,4)(3,6)\}$ & $3.579491$\\
\hline \hline
\multirow{4}{*}{\includegraphics[width=2.1cm]{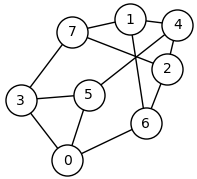} }& $\{(0,3,7,1,4,5)(2,6), (0,5,4,2,7,3)(1,6), (0,6)(1,7,2,4)(3,5)\}$ & $1.819247$\\
\cline{2-3}
& $\{(0,3)(1,6)(2,7)(4,5), (0,5,3,7,1,4,2,6), (0,6,2,4,1,7,3,5)\}$ & $1.983112$\\
\cline{2-3}
& $\{(0,3)(1,7,2,6)(4,5), (0,5)(1,6,2,4)(3,7), (0,6)(2,7,1,4)(3,5)\}$ & $2.418023$\\
\cline{2-3}
& $\{(0,3)(1,6)(2,7)(4,5), (0,5)(1,4)(2,6)(3,7), (0,6)(17)(2,4)(3,5)\}$ & $3.046722$\\
\hline\hline
\caption{$\tr(\widehat{M})$ of shunt-decomposition walks}
\end{longtable}

\section{Open Problems}
We have reformulated three models of discrete-time quantum walks in graph-theoretic language. As a first step to test the sensitivity of quantum walks to different combinatorial structures, we enumerated all rotation systems and shunt-decomposition of cubic graphs on up to $8$ vertices, and computed two parameters of the corresponding arc-reversal walks and shunt-decomposition walks. According to our data, we believe certain properties of the underlying graph structure are good indicators for the performance of potential quantum walk based algorithm. Conversely, one may also deduce useful information about the graph structure by running and measuring the associated quantum walk.

Here is a list of problems we would like to answer.
\begin{enumerate}
\item
When is the transition matrix $U$ walk-regular? When is the average mixing $\widehat{M}$ rational?
\item
There are many graph products, such as the Cartesian product, the tensor product and the strong product. Is there a relation between the quantum walks on two graphs and the quantum walk on their products, given properly chosen coins? 
\item
Using $3\times 3$ coins, we can implement arc-reversal quantum walks on the following chain of gems
\[\Phi(X), \Phi(\Phi(X))), \Phi(\Phi(\Phi(X))), \cdots\]
The sizes of these graphs increase quickly. Is there a way to assign the coins, such that the quantum walks on these graphs are determined by the quantum walk on $\Phi(X)$?

\item
For a regular graph $X$, is there a coin for which the entropy of the average mixing matrix distinguishes non-isomorphic rotation systems? 
\item
Which shunt-decompositions/embeddings give the maximum trace/entropy of the average mixing matrix?
\item 
Which shunt-decompositions/embeddings give the minimum hitting time/mixing time?
\item
Which shunt-decompositions/embeddings give perfect state transfer?
\item 
In the shunt-decomposition model, is there a relation between the cycle-structures of the shunts and the limiting distribution?
\item 
Can we extend the shunt-decomposition model to directed graphs where all vertices have the same out-degree?

\end{enumerate}

\section*{Acknowledgment}
C Godsil acknowledges support from NSERC grants RGPIN 50503-10325 and DAS 50503-10723.

\bibliographystyle{amsplain}
\bibliography{dqw.bib}

\end{document}